\documentclass{article}
\pdfoutput=1

\usepackage[margin=1in]{geometry}
\usepackage[utf8]{inputenc}
\usepackage[algoruled,linesnumbered]{algorithm2e}
\usepackage{amsfonts}       % blackboard math symbols
\usepackage{amsmath}
\usepackage{amssymb}
\usepackage{booktabs}       % professional-quality tables
\usepackage[font=small]{caption}
\usepackage{doi}
\usepackage{enumitem}
\usepackage{float}
\usepackage{graphicx}
\usepackage{hyperref}       % hyperlinks
\usepackage{microtype}      % microtypography
\usepackage{multirow}
\usepackage{nicefrac}       % compact symbols for 1/2, etc.
\usepackage{soul}
\usepackage{subfigure}
\usepackage{url}            % simple URL typesetting

\providecommand{\keywords}[1]
{
  \hspace{2.5em}\small
  \textbf{Keywords:} #1
}

\title{\textbf{Parallel QR Factorization of Block Low-Rank Matrices}}

\author{M. Ridwan Apriansyah$^1$ and Rio Yokota$^2$}
\date{%
    $^1$School of Computing, Tokyo Institute of Technology\\%
    \texttt{ridwan@rio.gsic.titech.ac.jp}\\%
    $^2$Global Scientific Information and Computing Center, Tokyo Institute of Technology\\%
    \texttt{rioyokota@gsic.titech.ac.jp}\\
}

\hypersetup{
pdftitle={Parallel QR Factorization of Block Low-Rank Matrices},
pdfsubject={},
pdfauthor={M. Ridwan Apriansyah and Rio Yokota},
pdfkeywords={Block Low-Rank matrix, QR factorization, Householder reflections, task-based execution},
}

\begin{document}
\maketitle

\begin{abstract}
We present two new algorithms for Householder QR factorization of Block Low-Rank (BLR) matrices: one that performs block-column-wise QR, and another that is based on tiled QR. We show how the block-column-wise algorithm exploits BLR structure to achieve arithmetic complexity of $\mathcal{O}(mn)$, while the tiled BLR-QR exhibits $\mathcal{O}(mn^{1.5})$ complexity. However, the tiled BLR-QR has finer task granularity that allows parallel task-based execution on shared memory systems. We compare the block-column-wise BLR-QR using fork-join parallelism with tiled BLR-QR using task-based parallelism. We also compare these two implementations of Householder BLR-QR with a block-column-wise Modified Gram-Schmidt (MGS) BLR-QR using fork-join parallelism, and a state-of-the-art vendor-optimized dense Householder QR in Intel MKL. For a matrix of size 131k $\times$ 65k, all BLR methods are more than an order of magnitude faster than the dense QR in MKL. Our methods are also robust to ill-conditioning and produce better orthogonal factors than the existing MGS-based method. On a CPU with 64 cores, our parallel tiled Householder and block-column-wise Householder algorithms show a speedup of 50 and 37 times, respectively.
\end{abstract}
\section{Introduction}
\label{sec:introduction}
QR factorization plays a central role in solving scientific and engineering problems. It factorizes a matrix $A \in \mathbb{R}^{m \times n}$ ($m \geq n$) into
\begin{equation}\label{eq:qr_fact}
  A = QR,
\end{equation}
where $Q \in \mathbb{R}^{m \times m}$ is orthogonal ($QQ^T=Q^TQ=I$) and $R \in \mathbb{R}^{m \times n}$ is upper triangular. 
QR factorization is well-known as a stable direct method to solve least squares problems \cite{Trefethen1997}.
It is also used as a stable solver for linear systems, possibly as an alternative to LU decomposition without pivoting \cite{Golub1996}, and has been used in efficient algorithms for computing polar decomposition \cite{Nakatsukasa2010_QWDH} and spectral decomposition \cite{NakatsukasaHigham2013}.

A wide range of problems in computational science requires factorizing dense matrices. Since traditional factorization methods require $\mathcal{O}(n^3)$, they have become the bottleneck for large-scale computation. Therefore, many techniques have been proposed to perform efficient factorization by exploiting the underlying structure of the matrices. A notable example is the low-rank structure arising from the discretization of integral equations, where the resulting full-rank matrices have been shown to possess many rank-deficient off-diagonal blocks \cite{Hackbusch1999}. This observation provides the basis of hierarchical low-rank representations \cite{Borm2003, Bebendorf2008, Aminfar2016_HODLR} that greatly reduce storage requirement and allow factorization in linear-polylogarithmic time. A similar low-rank structure has also been exploited in solving Toeplitz least squares problems \cite{Xi2014_ToeplitzLS} and sparse least squares problems \cite{Gnanasekaran2021_spaQR_LS, Gnanasekaran2021_spaQR_graphics}.

Block Low-Rank (BLR) matrices \cite{Amestoy2015} exploit a similar low-rank property, but produce a flat 2D blocked structure unlike the aforementioned hierarchical representations. This results in arithmetic complexity of $\mathcal{O}(n^2)$ for LU and QR factorization \cite{Amestoy2017,IdaNakashima2019}. Even though the hierarchical representations achieve lower complexity than BLR, the simplicity and flexibility of the BLR format make it easy to use in the context of a general-purpose, algebraic solver \cite{Amestoy2019}. Its simple, non-hierarchical structure is also efficient on parallel computers \cite{Jeannerod2019_BLR_LU,Akbudak2018,Chahara2018,Sergent2016}. Moreover, it has been shown that matrix-vector multiplication based on BLR-matrices is significantly faster than Hierarchical matrices for a large number of processes \cite{IdaNakashima2018}. For these reasons, BLR-matrices may be a better choice, at least for some problem classes and sizes \cite{Amestoy2019}.

In recent years, BLR factorization has generated considerable research interest. BLR-LU direct solvers have been shown to perform better than full-rank solvers in a sequential environment \cite{Amestoy2019}. Cholesky direct solvers using the BLR format have been utilized for weather modeling \cite{Akbudak2017} and PASTIX supernodal solver \cite{Pichon2017} on multicore architectures. BLR direct solvers have also been used in large-scale computation on distributed memory systems for seismic and electromagnetic \cite{Shantsev2017} and geospatial statistics problems \cite{Akbudak2018, Cao2020}.

However, QR factorization of BLR-matrices is not a well-studied problem. Ida et al. have combined the Gram Schmidt orthogonalization with BLR-matrix arithmetic to perform QR factorization on distributed memory systems \cite{IdaNakashima2019}. Even so, the method still relies on the traditional fork-join approach that has relatively large synchronization overhead. It may also suffer from numerical instability as it inherits the property of Gram-Schmidt iteration. In this article, we present two new algorithms for the QR factorization of BLR-matrices: one that performs block-column-wise QR based on the blocked Householder method \cite{Golub1996}, and another one that is based on the tiled QR \cite{Gunter2005}. Using the numerically stable Householder triangularization and BLR-matrix arithmetic, the block-column-wise algorithm achieves a theoretical complexity of $\mathcal{O}(mn)$, while the tiled BLR-QR exhibits $\mathcal{O}(mn^{1.5})$ complexity. Nonetheless, the tiled BLR-QR has finer granularity that allows for parallel task-based execution on shared memory systems. This leads to an out-of-order execution with very loose synchronization compared to the fork-join model. Numerical experiments show that our algorithms are more than an order of magnitude faster than the vendor-optimized dense Householder QR in Intel MKL. Moreover, our algorithms are robust to ill-conditioning and achieve higher parallel speedups compared to the existing Gram-Schmidt-based algorithm.

The rest of this article is organized as follows. In Section 2, we first summarize well-known methods to perform blocked and tiled QR decomposition of dense matrices. In Section 3 we briefly introduce BLR-matrices and elaborate on different methods to perform QR decomposition on them. We begin with an overview of the existing method that relies on modified Gram-Schmidt iteration. Then we explain our new algorithms: the first one that performs block-column-wise Householder QR; followed by one that is based on the tiled Householder QR. We then present the parallelization of our algorithms in Section 4, using both traditional fork-join and modern task-based execution models on shared-memory systems. Section 5 presents the results of various numerical experiments to show the performance and accuracy of our algorithms using several examples. Section 6 concludes this article.
\section{Block Dense QR}
\label{sec:block_dense_qr}
In this section, we summarize the standard blocked and tiled methods for QR decomposition of dense matrices. Although they have similar $O(mn^2)$ arithmetic complexity as the unblocked version, the blocked and tiled methods are known to be more efficient on modern supercomputers since they are rich in Level-3 BLAS operations that provide high performance on memory hierarchy systems \cite{Buttari2009_ParTiledDense}.

To facilitate this, we assume that the matrix $A \in \mathbb{R}^{m \times n}$ is subdivided into $p \times q$ square blocks of size $b \times b$, where $b$ is the chosen block size, $p=m/b$, and $q=n/b$. Extension to rectangular block is possible with extra permutation steps.
\subsection{Blocked Modified Gram-Schmidt Dense QR}
\label{subsec:dense_mbgs}
The modified Gram Schmidt (MGS) orthogonalization is a well-known method for computing QR factorization. The blocked version of this method is proposed in \cite{JalbyPhilippe1991_MBGS}. Consider a matrix $A \in \mathbb{R}^{m \times n} = [A_1 \;\; A_2]$ such that $A_1$ and $A_2$ are the block-columns of $A$. Then we can rewrite Equation \ref{eq:qr_fact} as
\begin{equation}
  \left[ A_1 \; \; A_2 \right] = \left[ Q_1 \; \; Q_2 \right]
  \left[ \begin{array}{cc}
    R_{1,1} & R_{1,2} \\
    \mathbf{0} & R_{2,2}
  \end{array} \right].
\end{equation}
Thus $A$ can be orthogonalized by the following steps:
\begin{enumerate}
  \item \textit{Orthogonalize block-column} $A_1 = Q_1 R_{1,1}$.
  \item \textit{Compute} $R_{1,2} = Q_{1}^{T} A_2$.
  \item \textit{Update} $A_2 \leftarrow A_2 - Q_1 R_{1,2}$.
  \item \textit{Orthogonalize block-column} $A_2 = Q_2 R_{2,2}$.
\end{enumerate}
These steps can be extended for arbitrary number of block columns, as shown in Algorithm \ref{alg:dense_mbgs_qr}. This is typically the method of choice when $A$ is well-conditioned. However, when $A$ is ill-conditioned, this method suffers from numerical instability due to rounding errors inherent in floating-point arithmetic on computers \cite{Trefethen1997}.
\begin{algorithm}[h]
  \DontPrintSemicolon
  \small
  \caption{Blocked Modified Gram-Schmidt (MGS) QR factorization}
  \label{alg:dense_mbgs_qr}
  \SetAlgoLined
  \KwIn{$A$ with $p \times q$ blocks}
  \KwOut{$Q$ with $p \times q$ blocks and $R$ with $q \times q$ blocks such that $A = QR$}
  \For{$j = 1$ \KwTo $q$} {
    $[Q_j, R_{j,j}]$ = QR($A_j$)\;
    \For{$k = j+1$ \KwTo $q$} {
        $R_{j,k} = Q_j^T A_{k}$\;
        $A_k \leftarrow A_k - Q_j R_{j,k}$\;
    }
  }
\end{algorithm}

\subsection{Blocked Householder Dense QR}
\label{subsec:dense_bch}
Householder triangularization is the principal method of QR factorization for its numerical stability. The blocked version of this method is proposed in \cite{Golub1996}, where the basic idea is to reorganize the computation by applying the Householder transformations in a cluster of columns at a time, i.e. triangularizing one block-column at a time. 

Before we proceed to the blocked version, let us recall the standard non-blocked method. Householder QR produces a factorization as in Equation \ref{eq:qr_fact} by performing Householder triangularization on $A$. However, unlike the Gram Schmidt QR, it does not directly produce the orthogonal factor $Q$. Instead, it performs in-place factorization such that in the end $A$ is replaced with $R$ in its upper triangular part and Householder vectors in its lower triangular part. These Householder vectors are then used to construct $Q$ when needed. This extra step of construction from Householder vectors is rich in matrix-vector multiplication and cannot fully utilize the parallelism in modern computers. For this reason, the compact WY representation \cite{Schreiber1989} has been proposed, which accumulates Householder reflectors such that
\begin{equation}\label{eq:compact_wy}
  Q=I-YTY^T,
\end{equation}
where $Y \in \mathbb{R}^{m \times n}$ is a unit lower trapezoidal matrix containing Householder vectors, and $T \in \mathbb{R}^{n \times n}$ is upper triangular. The generation of $Y$ and $T$ requires extra steps of $\mathcal{O}(n^3)$. Using this scheme, $Y$ can still be stored in the lower triangular part of $A$, but additional storage for $T$ is needed. Throughout this article, we assume that Householder QR factorization on dense matrices uses this representation to store $Q$.

For the sake of simplicity, consider the matrix $A \in \mathbb{R}^{m \times n}$ partitioned into a $3 \times 2$ block matrix
\begin{equation}\label{eq:a_blocked_2x2}
  A = 
  \left[ \begin{array}{cc}
    A_{1,1} & A_{1,2} \\
    A_{2,1} & A_{2,2} \\
    A_{3,1} & A_{3,2} \\
  \end{array} \right],
\end{equation}
where $A_{i,j} \in \mathbb{R}^{b \times b}$. If we rewrite Equation \ref{eq:qr_fact} as
\begin{equation}\label{eq:qr_blocked_2x2}
  \left[ \begin{array}{cc}
    A_{1,1} & A_{1,2} \\
    A_{2,1} & A_{2,2} \\
    A_{3,1} & A_{3,2}
  \end{array} \right] =
  Q \left[ \begin{array}{cc}
    R_{1,1} & R_{1,2} \\
    \mathbf{0} & R_{2,2} \\
    \mathbf{0} & \mathbf{0}
  \end{array} \right],
\end{equation}
$A$ can be triangularized as follows:
\begin{enumerate}
  \item \textit{Triangularize block-column}
  $\left[ \begin{array}{c}
    A_{1,1} \\
    A_{2,1} \\
    A_{3,1}
  \end{array} \right] = \hat{Q}_{1} \left[ \begin{array}{c}
    R_{1,1} \\
    \mathbf{0} \\
    \mathbf{0}
  \end{array} \right]$.
  \item \textit{Update}
  $\left[ \begin{array}{c}
    R_{1,2} \\
    A_{2,2} \\
    A_{3,2} 
  \end{array} \right] \leftarrow \hat{Q}_{1}^T \left[ \begin{array}{c}
    A_{1,2} \\
    A_{2,2} \\
    A_{3,2}
  \end{array} \right]$.
  \item \textit{Triangularize block-column}
  $\left[ \begin{array}{c}
    A_{2,2} \\
    A_{3,2}
  \end{array} \right] = \hat{Q}_{2} \left[ \begin{array}{c}
    R_{2,2} \\
    \mathbf{0}
  \end{array} \right]$.
\end{enumerate}
The orthogonal factor $Q$ can then be constructed as
\begin{equation}\label{eq:dense_bch_q_construction}
    Q = \hat{Q}_1
    \left[ \begin{array}{cc}
    I & \mathbf{0} \\
    \mathbf{0} & \hat{Q}_{2}
  \end{array} \right].
\end{equation}
Algorithm \ref{alg:dense_bch_qr} shows the generalization of blocked Householder QR based on step (1)-(3) above for an arbitrary number of blocks. Note that the Householder triangularization ultimately amounts to multiplication by $Q^T$ from the left. Thus, when needed, multiplication with $Q$ can be done by performing similar steps in the correct order and transposition. Algorithm \ref{alg:dense_bch_left_multiply_q} shows the method to multiply a matrix $C$ by $Q$ from the left, which has roughly the same cost as Algorithm \ref{alg:dense_bch_qr}.
\begin{algorithm}[ht]
  \DontPrintSemicolon
  \small
  \caption{Blocked Householder QR factorization}
  \label{alg:dense_bch_qr}
  \SetAlgoLined
  \KwIn{$A$ with $p \times q$ blocks}
  \KwOut{$Y,R$ with $p \times q$ blocks and $T$ with $1 \times q$ blocks such that $R$ is upper triangular and $Y,T$ contain intermediate orthogonal factors}
  \For{$k = 1$ \KwTo $q$} {
    QR
    $\left(\left[
        \begin{array}{c}
            A_{k,k} \\
            A_{k+1,k} \\
            \vdots \\
            A_{p,k}
        \end{array}
    \right]\right)$
    = $\hat{Q}_{k}
    \left[
        \begin{array}{c}
            R_{k,k} \\
            \mathbf{0} \\
            \vdots \\
            \mathbf{0}
        \end{array}
    \right]$, such that
    $\hat{Q}_{k} = I - 
    \left[
        \begin{array}{c}
            Y_{k,k} \\
            Y_{k+1,k} \\
            \vdots \\
            Y_{p,k}
        \end{array}
    \right] T_k
    \left[
        \begin{array}{c}
            Y_{k,k} \\
            Y_{k+1,k} \\
            \vdots \\
            Y_{p,k}
        \end{array}
    \right]^T$\;
    \For{$j = k+1$ \KwTo $q$} {
        Update $\left[
            \begin{array}{c}
                R_{k,j} \\
                A_{k+1,j} \\
                \vdots \\
                A_{p,j}
            \end{array}
        \right] \leftarrow \hat{Q}_{k}^{T}
        \left[
            \begin{array}{c}
                A_{k,j} \\
                A_{k+1,j} \\
                \vdots \\
                A_{p,j}
            \end{array}
        \right]$\;
    }
  }
\end{algorithm}
\begin{algorithm}[ht]
  \DontPrintSemicolon
  \small
  \caption{Left multiplication by Q from blocked Householder algorithm}
  \label{alg:dense_bch_left_multiply_q}
  \SetAlgoLined
  \KwIn{$Y$ with $p \times q$ blocks, T with $1 \times q$ blocks, $C$ with $p \times q$ blocks}
  \KwOut{$C \leftarrow QC$}
  \For{$k = q$ \textbf{downto} 1} {
    \For{$j = q$ \textbf{downto} $k$} {
        Update $\left[
            \begin{array}{c}
                C_{k,j} \\
                C_{k+1,j} \\
                \vdots \\
                C_{p,j}
            \end{array}
        \right] \leftarrow
        \left(
            I - 
            \left[
                \begin{array}{c}
                    Y_{k,k} \\
                    Y_{k+1,k} \\
                    \vdots \\
                    Y_{p,k}
                \end{array}
            \right] T_k
            \left[
                \begin{array}{c}
                    Y_{k,k} \\
                    Y_{k+1,k} \\
                    \vdots \\
                    Y_{p,k}
                \end{array}
            \right]^T
        \right)
        \left[
            \begin{array}{c}
                C_{k,j} \\
                C_{k+1,j} \\
                \vdots \\
                C_{p,j}
            \end{array}
        \right]$\;
    }
  }
\end{algorithm}

\subsection{Tiled Householder Dense QR}
\label{subsec:dense_bh}
Gunter and Van De Geijn decomposed the operations of blocked Householder QR further to obtain the tiled Householder method \cite{Gunter2005}. The method takes its root in the updating factorization technique \cite{Golub1996, Stewart1998}. From Equation \ref{eq:qr_blocked_2x2}, the upper triangularization of $A$ can also be done by:
\begin{enumerate}
  \item \textit{Upper triangularize} block $A_{1,1}= \hat{Q}_{1,1} R_{1,1}$.
  \item \textit{Compute} $R_{1,2} = \hat{Q}_{1,1}^{T} A_{1,2}$.
  \item \textit{Upper triangularize}
  $\left[ \begin{array}{c}
    R_{1,1} \\
    A_{2,1}
  \end{array} \right] = \hat{Q}_{2,1} \left[ \begin{array}{c}
    {R_{1,1}}^{\prime} \\
    \mathbf{0}
  \end{array} \right]$, which zeroes $A_{2,1}$ and updates $R_{1,1}$.
  \item \textit{Update} $\left[ \begin{array}{c}
    R_{1,2} \\
    A_{2,2}
  \end{array} \right] \leftarrow \hat{Q}_{2,1}^{T} \left[ \begin{array}{c}
    R_{1,2} \\
    A_{2,2}
  \end{array} \right]$.
  \item \textit{Upper triangularize}
  $\left[ \begin{array}{c}
    {R_{1,1}}^{\prime} \\
    A_{3,1}
  \end{array} \right] = \hat{Q}_{3,1} \left[ \begin{array}{c}
    {R_{1,1}}^{\prime\prime} \\
    \mathbf{0}
  \end{array} \right]$, which zeroes $A_{3,1}$ and updates ${R_{1,1}}^{\prime}$.
  \item \textit{Update} $\left[ \begin{array}{c}
    R_{1,2} \\
    A_{3,2}
  \end{array} \right] \leftarrow \hat{Q}_{3,1}^{T} \left[ \begin{array}{c}
    R_{1,2} \\
    A_{3,2}
  \end{array} \right]$.
  \item \textit{Upper triangularize} block $A_{2,2} = \hat{Q}_{2,2} R_{2,2}$.
  \item \textit{Upper triangularize}
  $\left[ \begin{array}{c}
    R_{2,2} \\
    A_{3,2}
  \end{array} \right] = \hat{Q}_{3,2} \left[ \begin{array}{c}
    {R_{2,2}}^{\prime} \\
    \mathbf{0}
  \end{array} \right]$, which zeroes $A_{3,2}$ and updates $R_{2,2}$.
\end{enumerate}

Algorithm \ref{alg:dense_bh_qr} shows the generalization of step (1)-(8) for an arbitrary number of blocks. Analogous to the blocked Householder method, the construction of $Q$ can also be done by performing similar steps in the correct order and transposition. Algorithm \ref{alg:dense_bh_left_multiply_q} shows the method to left-multiply by $Q$, which also has roughly the same cost as Algorithm \ref{alg:dense_bh_qr}.

It is important to see that the tiled Householder QR is similar to the blocked Householder QR, except that the upper triangularization of a block column is decomposed into multiple, smaller operations involving at most two blocks at a time. This improves granularity and data locality of the operations. However, this leads to a larger storage requirements as we need to store more "T" matrices from the intermediate QR factorizations. On the contrary, the same approach to decompose block-column operation is quite difficult to apply to the blocked MGS method since orthogonalization needs all information of the column.
\begin{algorithm}[h]
  \DontPrintSemicolon
  \small
  \caption{Tiled Householder QR factorization}
  \label{alg:dense_bh_qr}
  \SetAlgoLined
  \KwIn{$A$ with $p \times q$ blocks}
  \KwOut{$Y,T,R$ with $p \times q$ blocks such that $R$ is upper triangular and $Y,T$ contain intermediate orthogonal factors}
  \For{$k = 1$ \KwTo $q$} {
    QR 
    $(A_{k,k}) = \hat{Q}_{k,k} R_{k, k}$, such that
    $\hat{Q}_{k,k} = I - Y_{k,k} T_{k,k} Y_{k,k}^T$\;
    \For{$j = k+1$ \KwTo $q$} {
        $R_{k,j} = \hat{Q}_{k,k}^T A_{k,j}$\;
    }
    \For{$i = k+1$ \KwTo $p$} {
        QR
        $\left(\left[
            \begin{array}{c}
                R_{k,k} \\
                A_{i, k}
            \end{array}
        \right]\right) = \hat{Q}_{i,k} 
        \left[
            \begin{array}{c}
                R_{k,k}^{\prime} \\
                \mathbf{0}
            \end{array}
        \right]$, such that
        $\hat{Q}_{i,k} = I - 
        \left[
            \begin{array}{c}
                I \\
                Y_{i,k}
            \end{array}
        \right] T_{i,k}
        \left[
            \begin{array}{c}
                I \\
                Y_{i,k}
            \end{array}
        \right]^T$ \;
        \For{$j = k+1$ \KwTo $q$} {
            $\left[
                \begin{array}{c}
                    R_{k, j} \\
                    A_{i, j}
                \end{array}
            \right] \leftarrow \hat{Q}_{i,k}^T
            \left[
                \begin{array}{c}
                    R_{k, j} \\
                    A_{i, j}
                \end{array}
            \right]$\;
        }
    }
  }
\end{algorithm}
\begin{algorithm}[h]
  \DontPrintSemicolon
  \small
  \caption{Left multiplication by Q from tiled Householder algorithm}
  \label{alg:dense_bh_left_multiply_q}
  \SetAlgoLined
  \KwIn{$C,Y,T$ with $p \times q$ blocks}
  \KwOut{$C \leftarrow QC$}
  \For{$k = q$ \textbf{downto} $1$} {
    \For{$i = p$ \textbf{downto} $k+1$} {
        \For{$j = q$ \textbf{downto} $k$} {
            Update $\left[
                \begin{array}{c}
                    C_{k, j} \\
                    C_{i, j}
                \end{array}
            \right] \leftarrow
            \left(
                \left[
                    \begin{array}{c}
                        I \\
                        Y_{i,k}
                    \end{array}
                \right] T_{i,k}
                \left[
                    \begin{array}{c}
                        I \\
                        Y_{i,k}
                    \end{array}
                \right]^T
            \right)
            \left[
                \begin{array}{c}
                    C_{k, j} \\
                    C_{i, j}
                \end{array}
            \right]$\;
        }
    }
    \For{$j = q$ \textbf{downto} $k$} {
        Update $C_{k,j} \leftarrow (I - Y_{k,k} T_{k,k} Y_{k,k}^T) C_{k,j}$\;
    }
  }
\end{algorithm}
\section{Block Low-Rank QR}
\label{sec:blr_qr}
\subsection{Block Low-Rank Matrices}
\label{subsec:blr_matrices}
Block low-rank matrices exploit the low-rank property by performing flat 2D blocking and storing rank-deficient (\textit{admissible}) blocks in low-rank representation. We briefly introduce them in this section. For a detailed explanation, we refer the reader to \cite{Amestoy2015}.

Given a dense matrix $A \in \mathbb{R}^{m \times n}$, a block size $b$ is chosen to subdivide the matrix into $p \times q$ ($p=m/b$, $q=n/b$) square blocks such that
\begin{equation*}
  \tilde{A} = \left[
    \begin{array}{cccc}
        \tilde{A}_{1,1} & \tilde{A}_{1,2} & \cdots & \tilde{A}_{1,q} \\
        \tilde{A}_{2,1} & \cdots & \cdots & \vdots \\
        \vdots & \cdots & \cdots & \vdots \\
        \tilde{A}_{p,1} & \cdots & \cdots & \tilde{A}_{p,q}
    \end{array}
    \right],
\end{equation*}
where
\begin{equation}\label{eq:blr_matrix_def}
  \tilde{A}_{i,j} =
    \begin{cases}
      A_{i,j} & (A_{i,j} \text{ is not admissible}) \\
      U_{i,j} V_{i,j}^{T} & (A_{i,j} \text{ is admissible})
    \end{cases} \bigg| 1 \leq i \leq p,\,1 \leq j \leq q
\end{equation}
and $U_{i,j} V_{i,j}^{T}$ is the low-rank compressed form of $A_{i,j}$ so that $U_{i,j}, V_{i,j} \in \mathbb{R}^{b \times r_{i,j}}$.
In practice, the compression rank is chosen adaptively for each block. Given a prescribed error tolerance $\epsilon > 0$, we choose $r_{i,j} = r_{\epsilon}$ to be the smallest integer that satisfies $\|\tilde{A}_{i,j}-A_{i,j}\|_F \leq \epsilon \|A_{i,j}\|_F$, where $\| \cdot \|_F$ denotes the Frobenius norm. This implies the overall BLR compression error
\begin{equation}
    \|\tilde{A}-A\|_F \leq \epsilon \|A\|_F    
\end{equation}
is also bounded by $\epsilon$ (see \cite{Ida2015}).
Let us denote $r$ as the maximum rank of the admissible blocks. In many problem classes, $r$ can be shown to be much smaller than $n$ \cite{Bebendorf2003, Bebendorf2005}. Throughout this article we assume that $r$ is a small constant, i.e. $r=\mathcal{O}(1)$.

The admissibility condition determines whether a block is admissible for low-rank compression. There are mainly two types of admissibility conditions. One is the \textit{weak admissibility} where all off-diagonal blocks are assumed to be admissible. Another one is the \textit{strong admissibility} where inadmissible off-diagonal blocks exist. Under the strong admissibility condition, one typically use the admissibility constant $\eta > 0$ to determine admissible blocks based on the underlying geometric information of the matrix (see \cite{Amestoy2017} for details). Here we assume that the larger the $\eta$, the more inadmissible off-diagonal blocks exist.
When geometric information is not available, one approach is to attempt to compress each block and revert the ones with large rank (e.g. larger than $b/2$) back to dense form. This requires extra operations which is usually negligible in case of BLR compression.

\begin{figure}[h]
    \centering
    \subfigure{\includegraphics[width=0.30\linewidth]{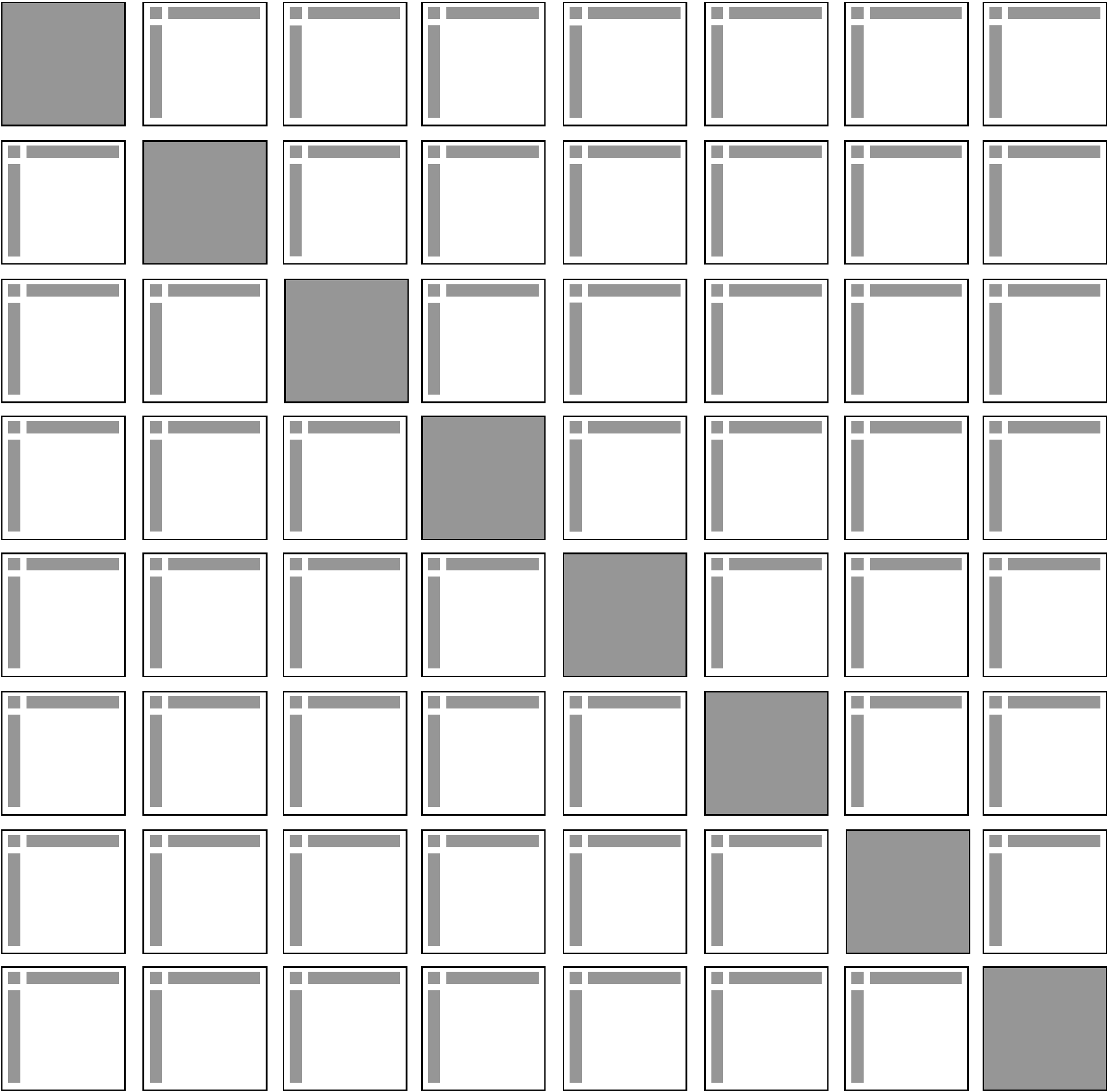}}%
    \hspace{0.05\linewidth}%
    \subfigure{\includegraphics[width=0.30\linewidth]{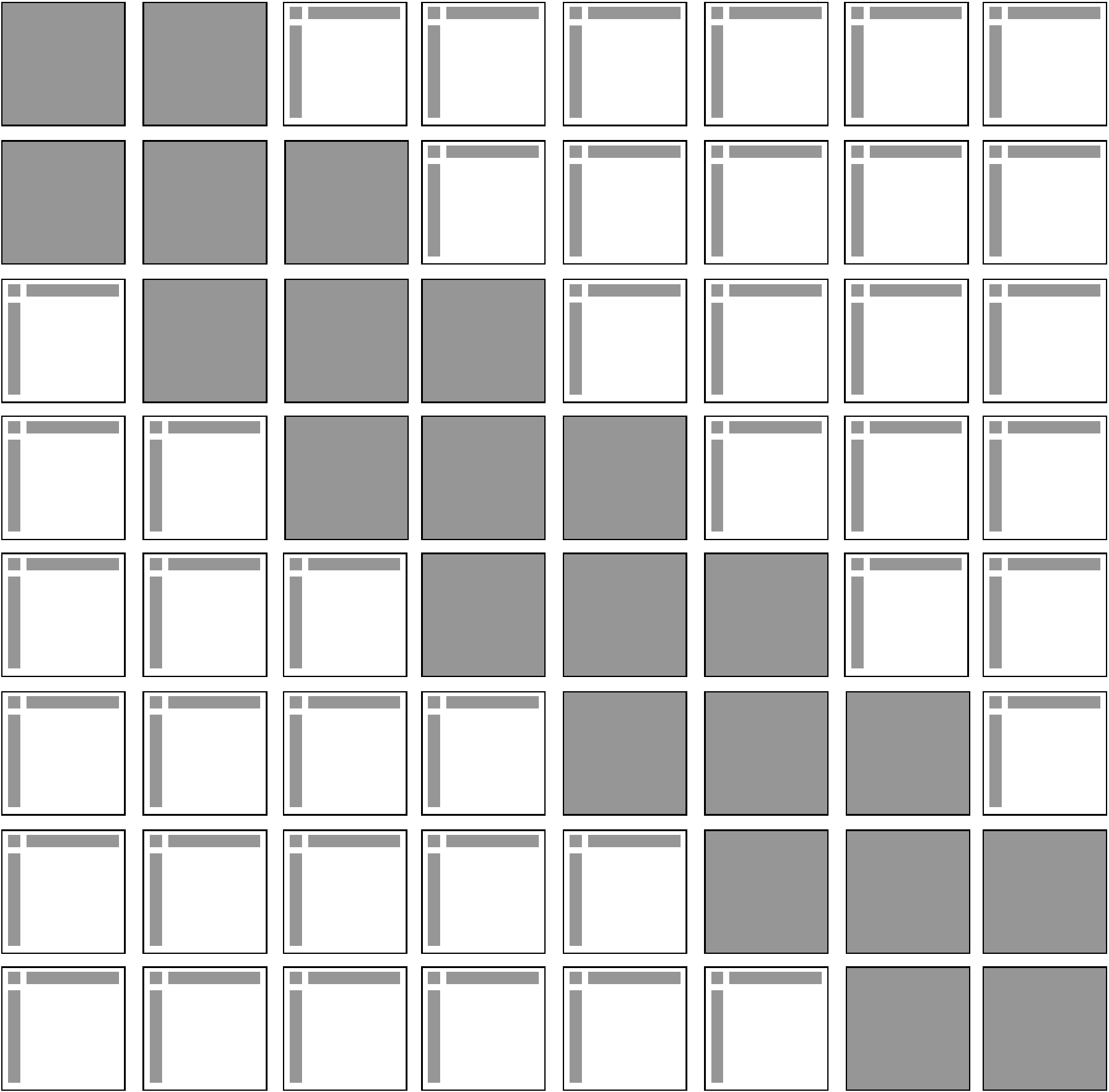}}
    \caption{Example of block low-rank matrices: weakly admissible (left); strongly admissible (right)}
    \label{fig:blr_matrices}
\end{figure}

The choice of admissibility condition depends on the problem. When all off-diagonal blocks of the matrix have small rank, e.g. in 2D Poisson problems, weak admissibility condition is usually sufficient. However when the matrix has full-rank off-diagonal blocks, e.g. in 3D Helmholtz problems, strong admissibility is commonly preferred. Figure \ref{fig:blr_matrices} shows the examples of BLR-matrices with different admissibility conditions. Grey-filled squares represent inadmissible (dense) blocks and the other represent admissible (low-rank) blocks.

The storage and computational costs involving BLR-matrices depend on the chosen block size and admissibility condition. Unlike hierarchical representations where the (minimum) block size only contributes to lower order terms of the cost, block size in BLR has a significant effect due to the flat blocking scheme. It has been shown in \cite{Amestoy2017} that the optimal choice of block size $b=\mathcal{O}(\sqrt{n})$ leads to $\mathcal{O}(n^{1.5})$ storage for a square BLR-matrix and $\mathcal{O}(n^2)$ cost for BLR factorization. Moreover, assuming that the cost to compress an admissible block is $\mathcal{O}(b^2r)$, which for example is achievable using Rank Revealing QR \cite{Golub1996} or randomized SVD \cite{HalkoMartinsson2011_RSVD}, constructing a BLR-matrix also requires $\mathcal{O}(n^2)$. Here we assume that each admissible block is compressed using Rank Revealing QR factorization. This produces the approximation $\tilde{A}_{i,j} = U_{i,j} V_{i,j}^T$ such that $U_{i,j}$ has orthonormal columns.

Operations on BLR-matrices are similar to that of blocked dense matrices, with the addition of basic operations involving dense and low-rank blocks. Let us define four operations:
\begin{itemize}[itemsep=0pt]
    \item[] {\ttfamily{OP1}}: Low-Rank $+$ Dense = Dense
    \item[] {\ttfamily{OP2}}: Low-Rank $+$ Low-Rank = Low-Rank
    \item[] {\ttfamily{OP3}}: Low-Rank $\times$ Dense = Low-Rank
    \item[] {\ttfamily{OP4}}: Low-Rank $\times$ Low-Rank = Low-Rank
\end{itemize}
{\ttfamily{OP1}} requires converting the low-rank block to dense form followed by dense blocks addition, which amounts to a total cost of $\mathcal{O}(b^2r)$.
In {\ttfamily{OP2}}, two low-rank blocks can be summed by simply agglomerating their components. However, this would make the rank grow very quickly. Therefore, the rounded addition method in Algorithm \ref{alg:lowrank_addition} \cite[p.~16]{Bebendorf2008} is employed to efficiently re-compress the resulting matrix back to rank $r_{\epsilon}$. Here we assume that recompression is performed every time two low-rank matrices are added, leading to a cost of $\mathcal{O}(br^2)$ per operation.
For {\ttfamily{OP3}}, the operation boils down to multiplication between the dense block and $V^T$ part of the low-rank block from the right (or $U$ part if it's from the left), which costs $\mathcal{O}(b^2r)$.
Lastly, for {\ttfamily{OP4}}, the resulting low-rank block $C=A \times B$ is formed by
\begin{displaymath}
    \underbrace{U_A}_{U_C} \underbrace{V_A^T U_B V_B^T}_{V_C^T},
\end{displaymath}
which costs $\mathcal{O}(br^2)$.
\begin{algorithm}[h]
  \DontPrintSemicolon
  \small
  \caption{Rounded Low-Rank Addition}
  \label{alg:lowrank_addition}
  \SetAlgoLined
  \KwIn{$A=U_AV_A^T$, $B=U_BV_B^T$, $U_A,V_A \in \mathbb{R}^{b \times r_A}$, $U_B,V_B \in \mathbb{R}^{b \times r_B}$}
  \KwOut{$A+B \approx C=U_CV_C^T$, $U_C,V_C \in \mathbb{R}^{b \times r_{\epsilon}}$}
  $U=\left[
    \begin{array}{cc}
         U_A & U_B
    \end{array}
  \right]$,
  $V=\left[
    \begin{array}{cc}
         V_A & V_B
    \end{array}
  \right]$\;
  $\left[ Q_U, R_U \right] = QR(U)$, $\left[ Q_V, R_V \right] = QR(V)$\;
  $\left[ u,\sigma,v \right] = SVD(R_U R_V^T)$\;
  $U_C \leftarrow$ first $r_{\epsilon}$ columns of $Q_U u$\;
  $V_C \leftarrow$ first $r_{\epsilon}$ columns of $Q_V v \sigma$ \;
\end{algorithm}
\subsection{Blocked Modified Gram-Schmidt BLR-QR}
\label{subsec:blr_mbgs}
Ida et al. \cite{IdaNakashima2019} have combined the blocked MGS method with BLR-matrix arithmetic to formulate a BLR-QR decomposition algorithm. Given a weakly admissible BLR-matrix $\tilde{A}$, the algorithm produces two BLR-matrices $\tilde{Q}$ and $\tilde{R}$ under two assumptions:
\begin{itemize}[leftmargin=*]
    \item $\tilde{Q}$ and $\tilde{R}$ are BLR-matrices with the same structure as $\tilde{A}$
    \item The off-diagonal blocks $\tilde{Q}_{i,j}$ and $\tilde{R}_{i,j}$ are approximated using the same error threshold as the corresponding $\tilde{A}_{i,j}$.
\end{itemize}
Under these conditions, the matrices $\tilde{Q}$ and $\tilde{R}$ serve as the approximate QR decomposition of $\tilde{A}$, such that $\tilde{A} \approx \tilde{Q}\tilde{R}$. 
The algorithm proceeds in the same way as Algorithm \ref{alg:dense_mbgs_qr} with some operations tailored for BLR-matrices. Line 2 of Algorithm \ref{alg:dense_mbgs_qr}, which corresponds to orthogonalization of a block column, is performed based on the method presented in \cite{Benner2010}. The matrix multiplications in line 4 and 5 are performed using BLR-matrix arithmetic.

The QR factorization of a block column is described as follows. First, let us write the block column $\tilde{A}_j$ as
\begin{equation}\label{eq:blr_mbgs_block_column}
    \tilde{A}_j = 
    \left[
        \begin{array}{c}
            \tilde{A}_{1,j} \\
            \tilde{A}_{2,j}\\
            \vdots \\
            \tilde{A}_{p,j}
        \end{array}
    \right] =
    \left[
        \begin{array}{c}
            \tilde{A}^U_{1,j} \tilde{A}^V_{1,j}\\
            \tilde{A}^U_{2,j} \tilde{A}^V_{2,j}\\
            \vdots \\
            \tilde{A}^U_{p,j} \tilde{A}^V_{p,j}
        \end{array}
    \right],
\end{equation}
where 
\begin{equation}\label{eq:blr_blockcolumn_left_orthogonal}
    \tilde{A}^U_{i,j} =
    \begin{cases}
      I_b & (\tilde{A}_{i,j} \text{ is not admissible}) \\
      U_{i,j} & (\tilde{A}_{i,j} \text{ is admissible})
    \end{cases},\:
    \tilde{A}^V_{i,j} =
    \begin{cases}
      \tilde{A}_{i,j} & (\tilde{A}_{i,j} \text{ is not admissible}) \\
      V_{i,j}^T & (\tilde{A}_{i,j} \text{ is admissible})
    \end{cases}
\end{equation}
for $i=1,2,\ldots,p$. Note that we have just written each $\tilde{A}_{i,j}$ in left-orthogonal form. Now since each $\tilde{A}^V_{i,j}$ is a dense block, we can write the matrix
\begin{equation}\label{eq:blr_mbgs_bj}
    B_j = 
    \left[
        \begin{array}{c}
            \tilde{A}^V_{1,j}\\
            \tilde{A}^V_{2,j}\\
            \vdots \\
            \tilde{A}^V_{p,j}
        \end{array}
    \right]
\end{equation}
and perform the dense MGS QR factorization $B_j = \tilde{Q}^B_{j} \tilde{R}^B_{j}$. Then we partition the matrix $\tilde{Q}^B_{j}$ back according to the subdivision of $B_j$ in Equation \ref{eq:blr_mbgs_bj}. Since we know that each $\tilde{A}^U_{i,j}$ has orthonormal columns, we can obtain the orthogonal factor as
\begin{equation}\label{eq:blr_mbgs_qj}
    \tilde{Q}_j =
    \left[
        \begin{array}{cccc}
            \tilde{A}^U_{1,j} & \mathbf{0} & \cdots & \mathbf{0} \\
            \mathbf{0} & \tilde{A}^U_{2,j} & \ddots & \mathbf{0} \\
            \vdots & \ddots & \ddots & \mathbf{0} \\
            \mathbf{0} & \cdots & \mathbf{0} & \tilde{A}^U_{p, j}
        \end{array}
    \right] \tilde{Q}^B_{j} 
    =
    \left[
        \begin{array}{c}
            \tilde{A}^U_{1,j} \tilde{Q}^B_{1,j}\\
            \tilde{A}^U_{2,j} \tilde{Q}^B_{2,j}\\
            \vdots \\
            \tilde{A}^U_{p,j} \tilde{Q}^B_{p,j}
        \end{array}
    \right],
\end{equation}
and upper triangular $\tilde{R}_{j,j} = \tilde{R}^B_{j}$ such that $\tilde{A}_j = \tilde{Q}_j \tilde{R}_{j,j}$, that is the QR factorization of the block column $\tilde{A}_j$.

Assuming a block size $b = \mathcal{O}(\sqrt{n})$, this algorithm has an arithmetic complexity of $\mathcal{O}(mn)$, which is faster than the $\mathcal{O}(mn^2)$ dense MGS algorithm. We refer the reader to \cite{IdaNakashima2019} for a detailed explanation.
\subsection{Blocked Householder BLR-QR}
\label{subsec:blr_bch}
Our first algorithm follows the blocked Householder dense QR in Section \ref{subsec:dense_bch} to perform orthogonal triangularization of BLR-matrices. However, a BLR-matrix may contain admissible (low-rank) blocks that need to be handled in a different way than inadmissible (dense) blocks. Thus, the operations need to be extended to handle these low-rank blocks. Let $\tilde{A}$ be a BLR-matrix as defined in Equation \ref{eq:blr_matrix_def}. Our algorithm uses the steps from Algorithm \ref{alg:dense_bch_qr} to produce the approximate QR factorization $\tilde{A} \approx \tilde{Q} \tilde{R}$ such that:
\begin{itemize}[leftmargin=*]
    \item $\tilde{Q}$ and $\tilde{R}$ are BLR-matrices with the same structure as $\tilde{A}$.
    \item The admissible blocks $\tilde{Q}_{i,j}$ and $\tilde{R}_{i,j}$ are approximated using the same error threshold as the corresponding admissible block $\tilde{A}_{i,j}$.
\end{itemize}

\subsubsection{Triangularization of block column}
\label{subsec:triangularize_block_column}
The first operation that we need to redefine is the QR factorization of a block-column (line 2 of Algorithm \ref{alg:dense_bch_qr}). We adopt a similar approach to the one presented in \cite{Kressner2018}. Let us write the $k$-th block column that needs to be triangularized as 
\begin{equation}\label{eq:blr_bch_blockcolumn}
    \left[
        \begin{array}{c}
            \tilde{A}_{k,k} \\
            \tilde{A}_{k+1,k} \\
            \vdots \\
            \tilde{A}_{p,k}
        \end{array}
    \right] =
    \left[
        \begin{array}{c}
            \tilde{A}^U_{k,k} \tilde{A}^V_{k,k}\\
            \tilde{A}^U_{k+1,k} \tilde{A}^V_{k+1,k}\\
            \vdots \\
            \tilde{A}^U_{p,k} \tilde{A}^V_{p,k}
        \end{array}
    \right]
\end{equation}
where $\tilde{A}^U_{i,k} \tilde{A}^V_{i,k}$ for $i=k,k+1,\ldots,p$ are the left-orthogonal forms of $\tilde{A}_{i,k}$ as defined in Equation \ref{eq:blr_blockcolumn_left_orthogonal}. Next we perform Householder triangularization to factorize the matrix
\begin{equation}\label{eq:blr_bch_col_qr}
    \left[
        \begin{array}{c}
            \tilde{A}^V_{k,k} \\
            \tilde{A}^V_{k+1,k} \\
            \vdots \\
            \tilde{A}^V_{p,k}
        \end{array}
    \right] = \tilde{Q}^V_{k} 
    \left[
        \begin{array}{c}
            \tilde{R}_{k,k} \\
            \mathbf{0} \\
            \vdots \\
            \mathbf{0}
        \end{array}
    \right],
\end{equation}
where 
\begin{equation}
    \tilde{Q}^V_{k} = I - 
    \left[
        \begin{array}{c}
            Y_{k,k} \\
            Y_{k+1,k} \\
            \vdots \\
            Y_{p,k}
        \end{array}
    \right] T_k
    \left[
        \begin{array}{c}
            Y_{k,k} \\
            Y_{k+1,k} \\
            \vdots \\
            Y_{p,k}
        \end{array}
    \right]^T
\end{equation}
such that $T_k \in \mathbb{R}^{b \times b}$ and $Y_{i,k}$ has the same dimension as $\tilde{A}^V_{i,k}$ for $i=k,k+1,\ldots,p$. Since each $\tilde{A}^U_{i,k}$ has orthonormal columns, we set $\tilde{Y}_{i,k} = \tilde{A}^U_{i,k} Y_{i,k}$ to obtain the orthogonal factor
\begin{equation}
    \hat{Q}_{k} = I - 
    \left[
        \begin{array}{c}
            \tilde{Y}_{k,k} \\
            \tilde{Y}_{k+1,k} \\
            \vdots \\
            \tilde{Y}_{p,k}
        \end{array}
    \right] T_k
    \left[
        \begin{array}{c}
            \tilde{Y}_{k,k} \\
            \tilde{Y}_{k+1,k} \\
            \vdots \\
            \tilde{Y}_{p,k}
        \end{array}
    \right]^T.
\end{equation}
It is important to note that $\tilde{Y}_{i,k}$ is a low-rank block if the corresponding $\tilde{A}_{i,k}$ is low-rank, otherwise it is dense.

To see why this works, let us define $W = diag(\tilde{A}^U_{k,k},\tilde{A}^U_{k+1,k},\allowbreak\ldots,\tilde{A}^U_{p,k})$. Because diagonal blocks are dense, $W=diag(I,\tilde{A}^U_{k+1,k},\allowbreak\ldots,\tilde{A}^U_{p,k})$. Multiplying $\hat{Q}_k^T$ to the block-column in Equation \ref{eq:blr_bch_blockcolumn} yields
\begin{align}\label{eq:blr_bch_blockcolumn_triangularization}
    \hat{Q}_{k}^T
    \left[
        \begin{array}{c}
            \tilde{A}_{k,k} \\
            \tilde{A}_{k+1,k} \\
            \vdots \\
            \tilde{A}_{p,k}
        \end{array}
    \right]
    &=
    \left(
        I - W
        \left[
            \begin{array}{c}
                {Y}_{k,k} \\
                {Y}_{k+1,k} \\
                \vdots \\
                {Y}_{p,k}
            \end{array}
        \right] T_k^T
        \left[
            \begin{array}{c}
                {Y}_{k,k} \\
                {Y}_{k+1,k} \\
                \vdots \\
                {Y}_{p,k}
            \end{array}
        \right]^T
        W^T
    \right)
    W
    \left[
        \begin{array}{c}
            \tilde{A}^V_{k,k} \\
            \tilde{A}^V_{k+1,k} \\
            \vdots \\
            \tilde{A}^V_{p,k}
        \end{array}
    \right] \nonumber \\
    &= W
    \left( \tilde{Q}^V_k \right)^T
    \left[
        \begin{array}{c}
            \tilde{A}^V_{k,k} \\
            \tilde{A}^V_{k+1,k} \\
            \vdots \\
            \tilde{A}^V_{p,k}
        \end{array}
    \right] = W
    \left[
        \begin{array}{c}
            \tilde{R}_{k,k} \\
            \mathbf{0} \\
            \vdots \\
            \mathbf{0}
        \end{array}
    \right] =
    \left[
        \begin{array}{c}
            \tilde{R}_{k,k} \\
            \mathbf{0} \\
            \vdots \\
            \mathbf{0}
        \end{array}
    \right].
\end{align}
Thus, we have upper triangularized the $k$-th block-column and at the same time obtained the QR factorization of it.

The cost for this operation is dominated by the QR factorization in Equation \ref{eq:blr_bch_col_qr}. Each $\tilde{A}^V_{i,k}$ has a small number of rows $r$ ($\ll b$) if $\tilde{A}_{i,k}$ is low-rank.
Under weak admissibility condition, this leads to dense QR factorization of a $(b+(p-k)r) \times b$ matrix that costs $\mathcal{O}(b^3+pb^2r)$. The cost is similar under strong admissibility condition as long as the number of dense blocks in a block-column is bounded by a constant.

\subsubsection{Apply block column reflector}
\label{subsec:apply_block_column_reflector}
This operation corresponds to line 4 of Algorithm \ref{alg:dense_bch_qr} where we multiply the resulting $\hat{Q}_k$ with a block-column. Let us write the operation as
\begin{align}
    \hat{Q}_{k}^{T}
    \left[
        \begin{array}{c}
            \tilde{A}_{k,j} \\
            \tilde{A}_{k+1,j} \\
            \vdots \\
            \tilde{A}_{p,j}
        \end{array}
    \right]
    &=
    \left[
        \begin{array}{c}
            \tilde{A}_{k,j} \\
            \tilde{A}_{k+1,j} \\
            \vdots \\
            \tilde{A}_{p,j}
        \end{array}
    \right] -
    \left(
        \left[
            \begin{array}{c}
                \tilde{Y}_{k,k} \\
                \tilde{Y}_{k+1,k} \\
                \vdots \\
                \tilde{Y}_{p,k}
            \end{array}
        \right] T_k^T
        \left[
            \begin{array}{c}
                \tilde{Y}_{k,k} \\
                \tilde{Y}_{k+1,k} \\
                \vdots \\
                \tilde{Y}_{p,k}
            \end{array}
        \right]^T
        \left[
            \begin{array}{c}
                \tilde{A}_{k,j} \\
                \tilde{A}_{k+1,j} \\
                \vdots \\
                \tilde{A}_{p,j}
            \end{array}
        \right]
    \right) \notag \\
    &=
    \left[
        \begin{array}{c}
            \tilde{A}_{k,j} \\
            \tilde{A}_{k+1,j} \\
            \vdots \\
            \tilde{A}_{p,j}
        \end{array}
    \right] -
    \left[
            \begin{array}{c}
                \tilde{Y}_{k,k} \\
                \tilde{Y}_{k+1,k} \\
                \vdots \\
                \tilde{Y}_{p,k}
            \end{array}
        \right] \left[T_k^T\right]
    \left[
        \tilde{Y}_{k,k}^T \tilde{A}_{k,j} + \tilde{Y}_{k+1,k}^T \tilde{A}_{k+1,j} + \cdots + \tilde{Y}_{p,k}^T \tilde{A}_{p,j}
    \right]. \label{eq:blr_bch_apply_bc_reflector}
\end{align}
The operation shown in Equation \ref{eq:blr_bch_apply_bc_reflector} consists of multiplication between dense and low-rank blocks, and accumulation by low-rank addition. Under the weak admissibility condition, this operation costs $\mathcal{O}(b^2r+ pbr^2)$. Under strong admissibility condition some block-columns require multiplication between dense blocks, leading to a cost of $\mathcal{O}(b^3 + b^2r + pbr^2)$.

\subsubsection{Algorithm and Cost Estimate}
\begin{figure}[h]
    \centering
    \includegraphics[width=0.8\textwidth]{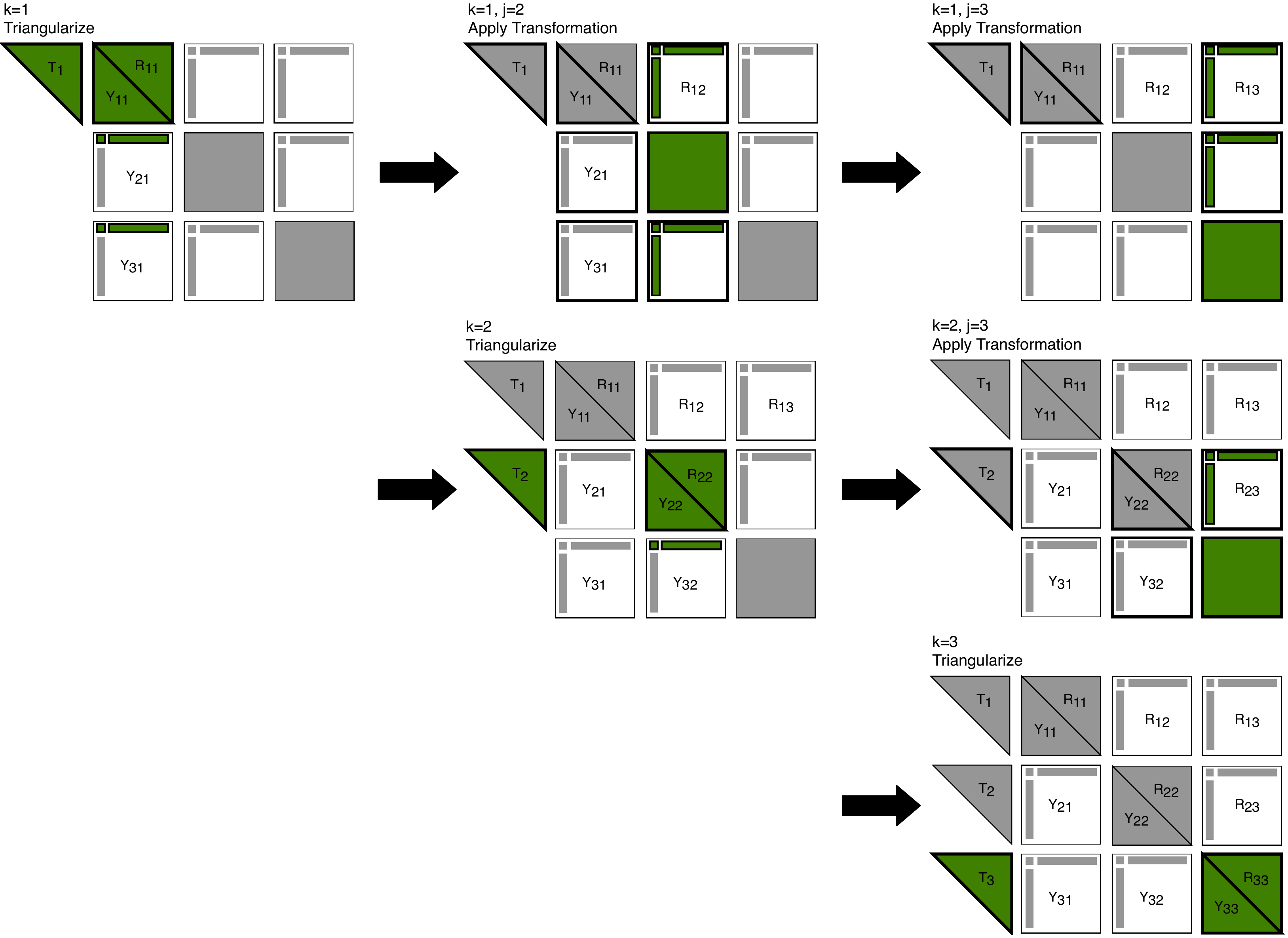}
    \caption{Graphical representation of operations in Algorithm \ref{alg:dense_bch_qr} on a BLR-matrix with $p=q=3$. Thick borders show the tiles that are being read and green fills shows the tiles that are being written at each step}
    \label{fig:blr_bch_illust}
\end{figure}

With the extended operations in place, we can perform blocked Householder QR to a BLR-matrix $\tilde{A}$. Figure \ref{fig:blr_bch_illust} shows the operations in Algorithm \ref{alg:dense_bch_qr} performed on a BLR-matrix with $3 \times 3$ blocks under weak admissibility condition.

In the following, we estimate the arithmetic complexity of our first algorithm. Under the strong admissibility condition, the structure of the BLR-matrix largely depends on the problem, requiring problem-specific analysis that is not the scope of this article. Thus, for simplicity let us assume the weak admissibility condition.

Table \ref{table:blr_bch_operations} shows the cost breakdown of one $k$-iteration in Algorithm \ref{alg:dense_bch_qr}. We get the total operation count by taking the sum for $k=1,2,\ldots,q$:
\begin{align*}
    T_{blocked}(m,n,b) &= \sum_{k=1}^{q} b^3 + pb^2r + (q-k)(b^2r + pbr^2) \notag\\
             &= \frac{1}{2} \left( 2nb^2 + 2mnr + n^2r + \frac{mn^2r^2}{b^2} - nbr - \frac{mnr^2}{b} \right).
\end{align*}
\begin{table}[h]
    \centering
    \caption{Operations inside one $k$-iteration ($1 \leq k \leq q$) of Algorithm \ref{alg:dense_bch_qr} and their costs}
    \label{table:blr_bch_operations}
    \small
    \begin{tabular}{lcc}
      \toprule
      \textbf{Operation} & \textbf{Complexity} & \textbf{Number of calls} \\
      \midrule
      Triangularization of block-column & $b^3+pb^2r$ & $1$ \\
      Apply block-column reflector & $b^2r+pbr^2$ & $q-k$ \\
      \bottomrule
    \end{tabular}
\end{table}

Setting $b = \mathcal{O}(\sqrt{n})$ yields the arithmetic complexity of $\mathcal{O}(mn)$.
In terms of storage, this algorithm uses the existing space of $\tilde{A}$ plus extra space to store $T_1, T_2, \ldots, T_q$ matrices, each of size $b \times b$. This extra storage is not larger than the space for a BLR-matrix, meaning that this algorithm requires $\mathcal{O}(m\sqrt{n})$ storage.
\subsection{Tiled Householder BLR-QR}
\label{subsec:blr_bh}
Our second algorithm is based on the tiled Householder dense QR explained in Section \ref{subsec:dense_bh}. It proceeds in the same way as Algorithm \ref{alg:dense_bh_qr} to produce an approximate QR factorization $\tilde{A} \approx \tilde{Q} \tilde{R}$ under the same assumption as our first algorithm. In the following, we explain how we extend the operations to handle low-rank blocks and estimate the overall cost of the algorithm.

\subsubsection{Householder QR factorization of diagonal block}
\label{subsec:geqrt}
This operation corresponds to line 2 of Algorithm \ref{alg:dense_bh_qr} where we perform Householder triangularization of a diagonal block
\begin{displaymath}
    \tilde{A}_{k,k} = \hat{Q}_{k,k} \tilde{R}_{k,k}.
\end{displaymath}
Since diagonal blocks are always dense, there is no need to handle low-rank blocks. This operation costs $\mathcal{O}(b^3)$.

\subsubsection{Apply block reflector}
\label{subsec:larfb}
This operation corresponds to line 4 of Algorithm \ref{alg:dense_bh_qr} where we multiply an off-diagonal block with orthogonal reflector, that is:
\begin{displaymath}
    \tilde{R}_{k,j} = \hat{Q}_{k,k}^T \tilde{A}_{k,j}.
\end{displaymath}
The off-diagonal block of a BLR-matrix can be dense or low-rank. If the particular block is low-rank, this operation amounts to a multiplication of dense and low-rank block that costs $\mathcal{O}(b^2r)$.

\subsubsection{Update QR factorization}
\label{subsec:tpqrt}
This operation corresponds to line 7 of Algorithm \ref{alg:dense_bh_qr} where we zero the off-diagonal block below $\tilde{A}_{k,k}$ using the QR factorization
\begin{displaymath}
    \text{QR}
    \left(\left[
        \begin{array}{c}
            \tilde{R}_{k,k} \\
            \tilde{A}_{i, k}
        \end{array}
    \right]\right) = \hat{Q}_{i,k} 
    \left[
        \begin{array}{c}
            \tilde{R}_{k,k}^{\prime} \\
            \mathbf{0}
        \end{array}
    \right].
\end{displaymath}
The upper triangular $\tilde{R}_{k,k}$ comes from the QR factorization of a diagonal block, so it is always a dense block. However $\tilde{A}_{i,k}$ can be a dense or low-rank block. In the case of a dense block, we simply concatenate the blocks and perform Householder dense QR on them. But if $\tilde{A}_{i,k}$ is low-rank, we first perform dense QR factorization of
\begin{equation}
    \left[
        \begin{array}{c}
            \tilde{R}_{k,k} \\
            \tilde{A}_{i,k}^V
        \end{array}
    \right] = \hat{Q}^V_{i,k} 
    \left[
        \begin{array}{c}
            \tilde{R}_{k,k}^{\prime} \\
            \mathbf{0}
        \end{array}
    \right],
\end{equation}
where
\begin{equation}
    \tilde{Q}^V_{i,k} = I - 
    \left[
        \begin{array}{c}
            I \\
            Y_{i,k}
        \end{array}
    \right] T_{i,k}
    \left[
        \begin{array}{c}
            I \\
            Y_{i,k}
        \end{array}
    \right]^T.
\end{equation}
Because $U_{i,k}$ has orthonormal columns, we set $\tilde{Y}_{i,k} = U_{i,k} Y_{i,k}$ to obtain the orthogonal factor
\begin{equation}\label{eq:trapezoidal_block_reflector}
    \hat{Q}_{i,k} = I - 
    \left[
        \begin{array}{c}
            I \\
            \tilde{Y}_{i,k}
        \end{array}
    \right] T_{i,k}
    \left[
        \begin{array}{c}
            I \\
            \tilde{Y}_{i,k}
        \end{array}
    \right]^T.
\end{equation}
Note that this is a specialization of the operation explained in Section \ref{subsec:triangularize_block_column} where the block column is composed of a dense upper triangular block on top of an off-diagonal block. This operation requires $\mathcal{O}(b^3)$ for both dense and low-rank $\tilde{A}_{i,k}$ due to the cost for generating $T_{i,k}$.

\subsubsection{Apply trapezoidal block reflector}
\label{subsec:tpmqrt}
This operation corresponds to line 9 of Algorithm \ref{alg:dense_bh_qr} where we multiply the orthogonal reflector from Equation \ref{eq:trapezoidal_block_reflector} to the corresponding blocks, that is
\begin{align*}
    \left[
        \begin{array}{c}
            \tilde{R}_{k, j} \\
            \tilde{A}_{i, j}
        \end{array}
    \right] \leftarrow \hat{Q}_{i,k}^T
    \left[
        \begin{array}{c}
            \tilde{R}_{k, j} \\
            \tilde{A}_{i, j}
        \end{array}
    \right]
    &= 
    \left[
        \begin{array}{c}
            \tilde{R}_{k, j} \\
            \tilde{A}_{i, j}
        \end{array}
    \right] -
    \left(
        \left[
            \begin{array}{c}
                I \\
                \tilde{Y}_{i,k}
            \end{array}
        \right] T_{i,k}^T
        \left[
            \begin{array}{c}
                I \\
                \tilde{Y}_{i,k}
            \end{array}
        \right]^T
        \left[
            \begin{array}{c}
                \tilde{R}_{k, j} \\
                \tilde{A}_{i, j}
            \end{array}
        \right]
    \right) \notag \\
    &=
    \left[
        \begin{array}{c}
            \tilde{R}_{k, j} \\
            \tilde{A}_{i, j}
        \end{array}
    \right] -
    \left[
        \begin{array}{c}
            I \\
            \tilde{Y}_{i,k}
        \end{array}
    \right] \left[ T_{i,k}^T \right]
    \left[
        \tilde{R}_{k,j} + \tilde{Y}_{i,k}^T \tilde{A}_{i,j}
    \right].
\end{align*}
The cost of this operation depends on the blocks $\tilde{Y}_{i,k}$, $\tilde{R}_{k,j}$, and $\tilde{A}_{i,j}$. If $\tilde{R}_{k,j}$ is low-rank and at least one of $\{ \tilde{Y}_{i,k}, \tilde{A}_{i,j} \}$ is low-rank, the cost is $\mathcal{O}(b^2r)$; Otherwise it is $\mathcal{O}(b^3)$.

\subsubsection{Algorithm and Cost Estimate}
\begin{figure}[h]
    \centering
    \includegraphics[width=0.8\textwidth]{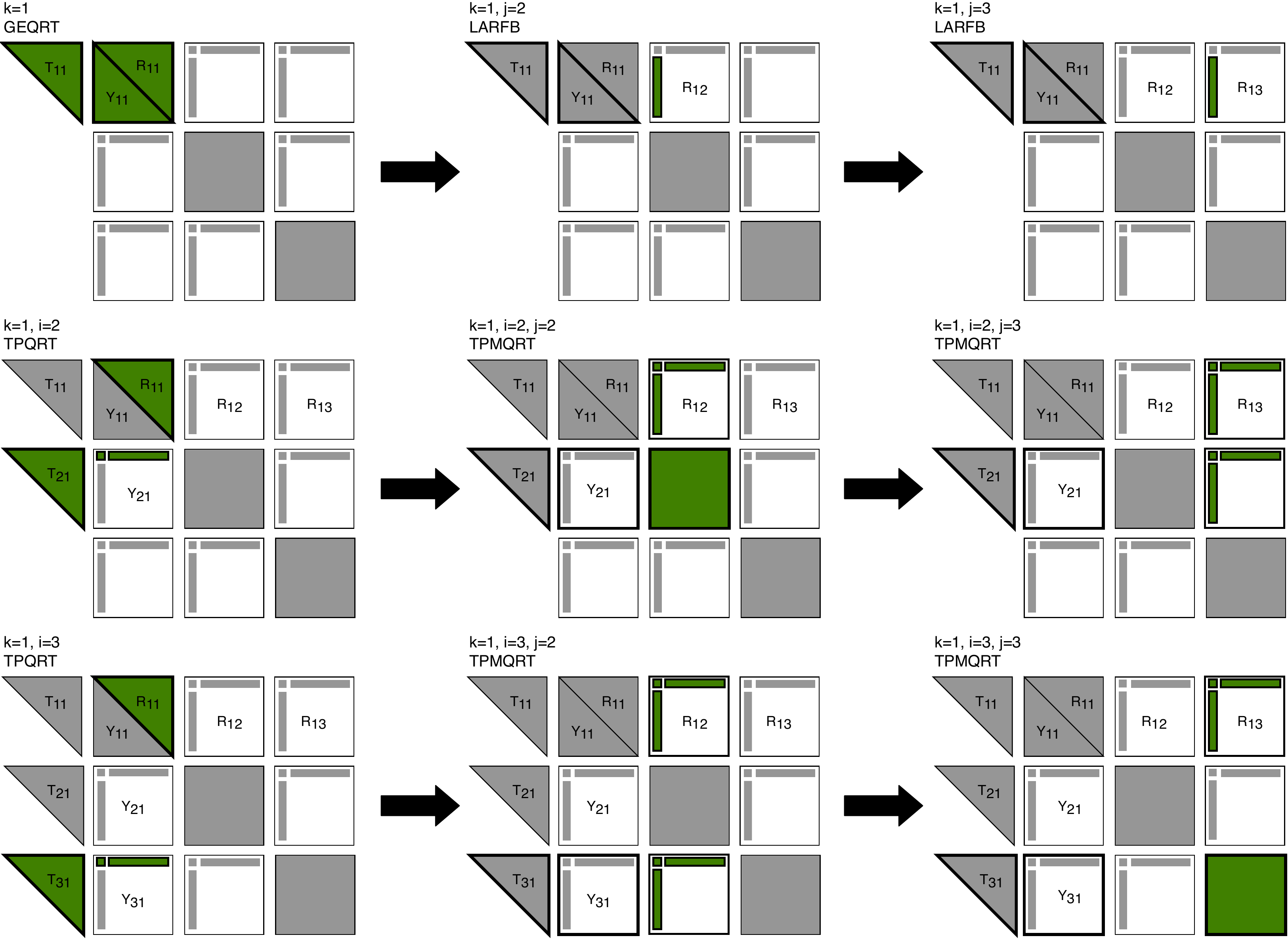}
    \caption{Graphical representation of one iteration of the outer loop in Algorithm \ref{alg:dense_bh_qr} on a BLR-matrix with $p=q=3$. Thick borders show the tiles that are being read and green fills shows the tiles that are being written at each step}
    \label{fig:blr_bh_illust}
\end{figure}

Once we have the extended operations defined, we can perform tiled Householder QR on a BLR-matrix $\tilde{A}$. Figure \ref{fig:blr_bh_illust} shows the operations that happen inside one outer iteration of Algorithm \ref{alg:dense_bh_qr} on a BLR-matrix with $3 \times 3$ blocks under weak admissibility condition.

In the following, we estimate the arithmetic complexity of our second algorithm. Let us also assume the weak admissibility condition for the sake of simplicity. Table \ref{table:blr_bh_operations} shows the cost breakdown of one $k$-iteration of Algorithm \ref{alg:dense_bh_qr}. Summing up for $k=1,2,\ldots,q$ gets us the total operation count:
\begin{align*}
    T_{tiled}(m,n,b) &= \sum_{k=1}^{q} (p-k+1)b^3 + (p-k+1)(q-k)b^2r \notag\\
             &= \frac{1}{6} \left( 6mnb + 3nb^2 + \frac{3mn^2r}{b} - 3n^2b - \frac{5n^3r}{b} - 3mnr  \right)
\end{align*}
\begin{table}[h]
    \centering
    \caption{Operations inside one $k$-iteration ($1 \leq k \leq q$) of Algorithm \ref{alg:dense_bh_qr} and their costs}
    \label{table:blr_bh_operations}
    \small
    \begin{tabular}{lcc}
      \toprule
      \textbf{Operation} & \textbf{Complexity} & \textbf{Number of calls} \\
      \midrule
      QR factorization of diagonal block & $b^3$ & $1$ \\
      Apply block reflector & $b^2r $ & $q - k$ \\
      Update QR factorization & $b^3 $ & $p - k$ \\
      Apply trapezoidal block reflector & $b^2r $ & $(p-k)(q-k)$ \\
      \bottomrule
    \end{tabular}
\end{table}

For $b = \mathcal{O}(\sqrt{n})$, this algorithm has an arithmetic complexity of $\mathcal{O}(mn^{1.5})$, which is slower than the blocked Householder variant. It also produces more $T$ matrices, which in total amounts to storing a lower trapezoidal $m \times n$ matrix. This leads to a storage requirement that grows similarly to that of dense Householder factorization ($\mathcal{O}(mn)$). However, this algorithm has finer granularity that makes it more efficient for parallel computation, which will be shown in later section.
\section{Multithreaded Block Low-Rank QR}
\label{sec:multithreaded_blr_qr}
In order to fully utilize modern multi-core architectures, we present the parallel algorithms for the BLR-QR factorization on shared memory systems. We first recall the fork-join parallelization of the blocked MGS-based method \cite{IdaNakashima2019}. Then we use a similar fork-join approach to parallelize our proposed algorithms. Lastly, we show the task-based parallel version of the tiled Householder BLR-QR, which is the main advantage of this variant.

Note that it is also possible to use task-based execution for the blocked MGS and Householder-based methods. The approach described in \cite{Kurzak2006} consists of representing the algorithm as a Directed Acyclic Graph (DAG) where nodes represent tasks (either block-column factorization or update) and edges represent dependencies among them. They refer this as dynamic lookahead technique. However, results show that their technique is still exposed to scalability problems due to the relatively coarse granularity of the tasks. Therefore we focus on using task-based execution on the finely grained tiled Householder-based algorithm.
\subsection{Parallel Blocked MGS BLR-QR}
\label{subsec:forkjoin_blr_mbgs}
Algorithm \ref{alg:dense_mbgs_qr} can be executed in parallel using the fork-join model \cite{Conway1963}, which has been presented in \cite{IdaNakashima2019}.
The block-column QR (line 2) can utilize a multithreaded BLAS/LAPACK kernel.
The computation of $R_{j,k}$ for different $k$ (line 4) can be performed in parallel. After that, the updates to $A_{k}$ can be performed in parallel too. These updates are first decomposed into independent block-by-block multiplications and computed simultaneously. Algorithm \ref{alg:forkjoin_blr_mbgs_qr} shows the fork-join parallel version of Algorithm \ref{alg:dense_mbgs_qr}.
\begin{algorithm}[h]
  \DontPrintSemicolon
  \small
  \caption{Fork-join blocked MGS BLR-QR factorization}
  \label{alg:forkjoin_blr_mbgs_qr}
  \SetAlgoLined
  \KwIn{$\tilde{A}$ with $p \times q$ blocks}
  \KwOut{$\tilde{Q}$ with $p \times q$ blocks and $\tilde{R}$ with $q \times q$ blocks such that $\tilde{A} \approx \tilde{Q}\tilde{R}$}
  \For{$j = 1$ \KwTo $q$} {
    $[\tilde{Q}_j, \tilde{R}_{j,j}]$ = QR($\tilde{A}_j$) \tcp{multithreaded}
    \For(\textbf{in parallel}){$k = j+1$ \KwTo $q$} {
        $\tilde{R}_{j,k} = \tilde{Q}_j^T \tilde{A}_{k}$\;
    }
    \ForAll(\textbf{in parallel}){ $\tilde{A}_{i,k}$ where $k > j$ and $1 \leq i
      \leq p$} {
      $\tilde{A}_{i,k} \leftarrow \tilde{A}_{i,k} - \tilde{Q}_{i,k}\tilde{R}_{j,k}$\;
    }
  }
\end{algorithm}
\subsection{Parallel Blocked Householder BLR-QR}
\label{subsec:parallel_blr_bch}
The fork-join approach can also be used to parallelize Algorithm \ref{alg:dense_bch_qr}. Let us look at the dependency among the operations. Operations in line 4 for different $j$ update different block-columns and only have a common dependency to the computation of $\hat{Q}_{k}$ (line 2). Thus they can be computed simultaneously as soon as the computation of $\hat{Q}_{k}$ is done.
The computation of $\hat{Q}_{k}$ itself (line 2) can be done by utilizing a multithreaded BLAS/LAPACK kernel.
Algorithm \ref{alg:forkjoin_blr_bch} shows the parallel version where the $j$-loop (line 3) branches off to become a parallel region. Similarly, the left multiplication by $\tilde{Q}$ in Algorithm \ref{alg:dense_bch_left_multiply_q} can also be executed in parallel using the same approach.
\begin{algorithm}[h]
  \DontPrintSemicolon
  \small
  \caption{Fork-join blocked Householder BLR-QR factorization}
  \label{alg:forkjoin_blr_bch}
  \SetAlgoLined
  \KwIn{$\tilde{A}$ with $p \times q$ blocks}
  \KwOut{$\tilde{Y},\tilde{R}$ with $p \times q$ blocks and $T$ with $1 \times q$ blocks such that $R$ is upper triangular and $\tilde{Y},T$ contain intermediate orthogonal factors}
  \For{$k = 1$ \KwTo $q$} {
    QR
    $\left(\left[
        \begin{array}{c}
            \tilde{A}_{k,k} \\
            \tilde{A}_{k+1,k} \\
            \vdots \\
            \tilde{A}_{p,k}
        \end{array}
    \right]\right)$
    = $\hat{Q}_{k}
    \left[
        \begin{array}{c}
            \tilde{R}_{k,k} \\
            \mathbf{0} \\
            \vdots \\
            \mathbf{0}
        \end{array}
    \right]$, such that
    $\hat{Q}_{k} = I - 
    \left[
        \begin{array}{c}
            \tilde{Y}_{k,k} \\
            \tilde{Y}_{k+1,k} \\
            \vdots \\
            \tilde{Y}_{p,k}
        \end{array}
    \right] T_k
    \left[
        \begin{array}{c}
            \tilde{Y}_{k,k} \\
            \tilde{Y}_{k+1,k} \\
            \vdots \\
            \tilde{Y}_{p,k}
        \end{array}
    \right]^T$ \tcp{\footnotesize multithreaded}
    \For(\textbf{in parallel}){$j = k+1$ \KwTo $q$} {
        Update $\left[
            \begin{array}{c}
                \tilde{R}_{k,j} \\
                \tilde{A}_{k+1,j} \\
                \vdots \\
                \tilde{A}_{p,j}
            \end{array}
        \right] \leftarrow \hat{Q}_{k}^{T}
        \left[
            \begin{array}{c}
                \tilde{A}_{k,j} \\
                \tilde{A}_{k+1,j} \\
                \vdots \\
                \tilde{A}_{p,j}
            \end{array}
        \right]$\;
    }
  }
\end{algorithm}
\subsection{Parallel Tiled Householder BLR-QR}
\label{subsec:parallel_blr_bh}
Looking at the operations in Algorithm \ref{alg:dense_bh_qr}, the fork-join approach is also applicable to obtain a parallel algorithm. The computations of $R_{k,j}$ in line 4 for different $j$ only have a common dependency to the computation of $\hat{Q}_{k,k}$ in line 2, and thus can be performed in parallel. A similar dependency can be seen among the update operations in line 9. Algorithm \ref{alg:forkjoin_blr_bh} shows the fork-join parallel version of Algorithm \ref{alg:dense_bh_qr} where the two $j$-loops (line 3 and 8) branch off to become parallel regions. 
Lastly, the computation of $\hat{Q}_{k,k}$ (line 2) and $\hat{Q}_{i,k}$ (line 7) can be performed using a multithreaded BLAS/LAPACK kernel.

\begin{algorithm}[h]
  \DontPrintSemicolon
  \small
  \caption{Fork-join tiled Householder BLR-QR factorization}
  \label{alg:forkjoin_blr_bh}
  \SetAlgoLined
  \KwIn{$\tilde{A}$ with $p \times q$ blocks}
  \KwOut{$\tilde{Y}, T, \tilde{R}$ with $p \times q$ blocks such that $\tilde{R}$ is upper triangular and $\tilde{Y},T$ contain intermediate orthogonal factors}
  \For{$k = 1$ \KwTo $q$} {
    QR
    $(\tilde{A}_{k,k}) = \hat{Q}_{k,k} \tilde{R}_{k, k}$, such that
    $\hat{Q}_{k,k} = I - \tilde{Y}_{k,k} T_{k,k} \tilde{Y}_{k,k}^T$ \tcp{multithreaded}
    \For(\textbf{in parallel}){$j = k+1$ \KwTo $q$} {
        $\tilde{R}_{k,j} = \hat{Q}_{k,k}^T \tilde{A}_{k,j}$\;
    }
    \For{$i = k+1$ \KwTo $p$} {
        QR
        $\left(\left[
            \begin{array}{c}
                \tilde{R}_{k,k} \\
                \tilde{A}_{i, k}
            \end{array}
        \right]\right) = \hat{Q}_{i,k} 
        \left[
            \begin{array}{c}
                \tilde{R}_{k,k}^{\prime} \\
                \mathbf{0}
            \end{array}
        \right]$, such that
        $\hat{Q}_{i,k} = I - 
        \left[
            \begin{array}{c}
                I \\
                \tilde{Y}_{i,k}
            \end{array}
        \right] T_{i,k}
        \left[
            \begin{array}{c}
                I \\
                \tilde{Y}_{i,k}
            \end{array}
        \right]^T$ \tcp{multithreaded}
        \For(\textbf{in parallel}){$j = k+1$ \KwTo $q$} {
            $\left[
                \begin{array}{c}
                    \tilde{R}_{k, j} \\
                    \tilde{A}_{i, j}
                \end{array}
            \right] \leftarrow \hat{Q}_{i,k}^T
            \left[
                \begin{array}{c}
                    \tilde{R}_{k, j} \\
                    \tilde{A}_{i, j}
                \end{array}
            \right]$\;
        }
    }
  }
\end{algorithm}

As we have mentioned before, the tiled Householder QR has finer granularity compared to the other blocked algorithms, which could be leveraged in a parallel environment. The idea of exploiting finer granularity of tiled Householder QR to obtain a highly parallel algorithm has been introduced in \cite{Buttari2009_ParTiledDense} in the context of optimizing dense factorization. Since we are using the same tiled algorithm but adapted to the BLR format, we follow similar steps to reach an efficient parallelization scheme of our BLR-QR algorithm.

\begin{figure}[h]
    \centering
    \includegraphics[width=0.6\textwidth]{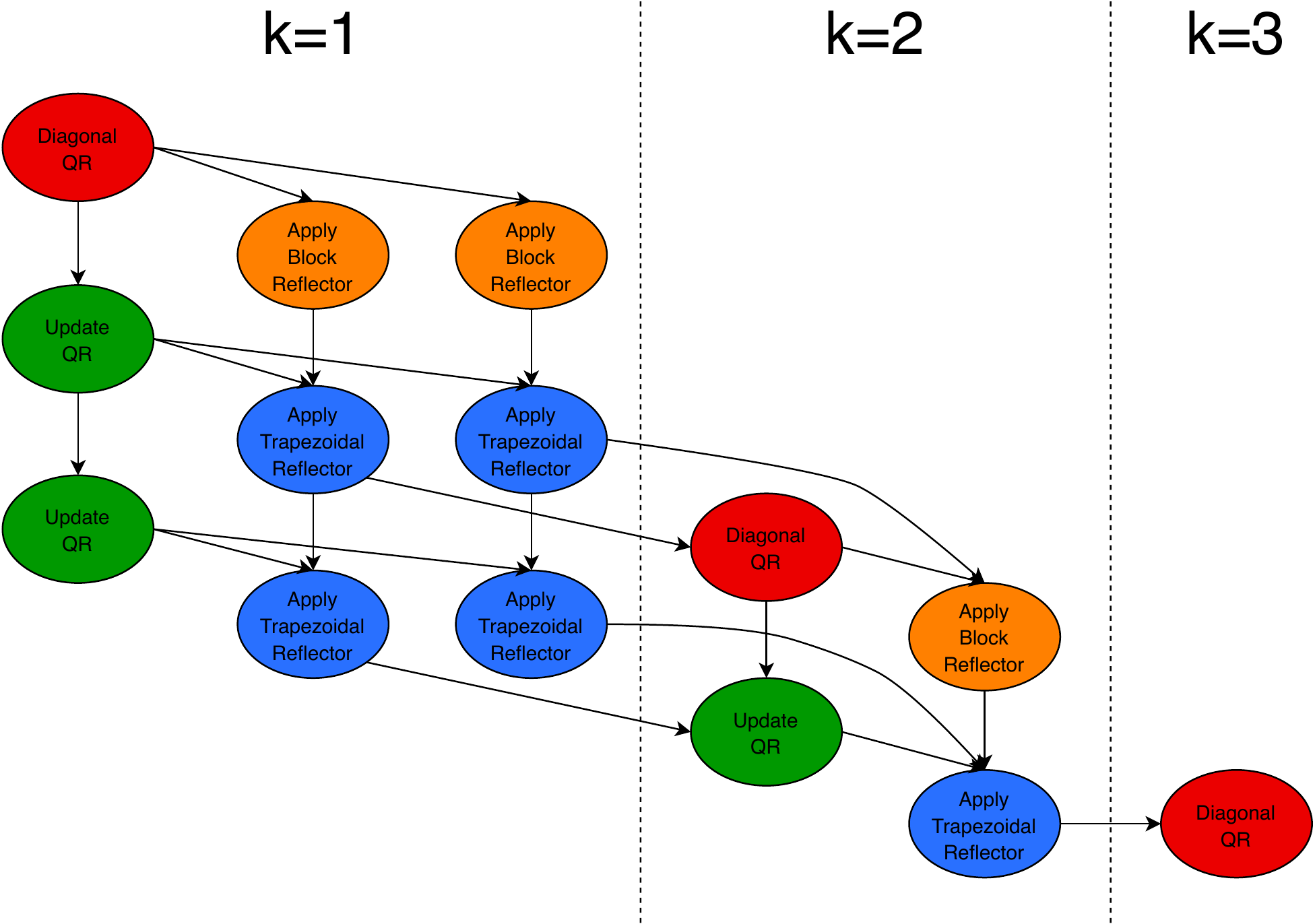}
    \caption{Dependency graph of Algorithm \ref{alg:dense_bh_qr} on a BLR-matrix with $p=q=3$}
    \label{fig:blr_bh_dag}
\end{figure}

Let us again look at the operations in Algorithm \ref{alg:dense_bh_qr}, but this time without limiting ourselves to one $k$ iteration. It turns out that there are direct dependencies between operations across consecutive iterations of $k$. This chain of dependencies is best described by a DAG, where a node represents an operation and an edge represents a dependency between two operations. Figure \ref{fig:blr_bh_dag} shows the dependency graph when Algorithm \ref{alg:dense_bh_qr} is executed on a BLR-matrix with $p=q=3$. Note that the construction of $\tilde{Q}$ uses the same set of operations so we can obtain a similar dependency graph from it.

It can be seen from Figure \ref{fig:blr_bh_dag} that the DAG also has a recursive structure. For any $p_1 \geq p_2$, the DAG for BLR-matrix with $p_2 \times p_2$ blocks is a subgraph of the DAG for BLR-matrix with $p_1 \times p_1$ blocks. This property allows for reusing the existing DAG to accelerate the construction of a larger graph.

Once we obtain the DAG, we can use it as a guide in executing the tasks. A task can be started as soon as all of its dependencies are fulfilled. Once a task $T$ is finished, the scheduler fulfills the dependency of tasks that are dependent on $T$. As a result, the threads only need to check the task pool and execute tasks that are "ready" to be executed. All threads repeat this cycle until the task pool is empty and the algorithm is finished. This results in an out-of-order execution with very loose synchronization required between the threads compared to the fork-join model.

A closer look into the graph in Figure \ref{fig:blr_bh_dag} reveals that certain kinds of task have more outgoing edges than the others. This offers a chance for improvement because executing a task with a large number of outgoing edges will fulfill more dependencies, thus bringing more tasks into the ``ready'' state. Therefore, a priority value can be assigned to each task, such that a task with more outgoing edges has a higher priority to be executed first. In our algorithm, it is clear that the factorization of diagonal blocks has the highest priority. The second one is updating QR factorization, followed by the application of trapezoidal and block reflectors.
\section{Numerical Results}
\label{sec:numerical_result}
In this section, we demonstrate the performance and accuracy of our algorithms using several example matrices on a shared-memory system. The BLR-QR algorithms were implemented in C++ where floating point calculations were performed in double precision. Fork-join and task-based parallelism were performed using OpenMP.
BLAS and LAPACK routines from Intel MKL were used for the inner kernels involving dense matrices, where single-threaded kernels were used inside OpenMP parallel regions and multi-threaded kernels were used outside of parallel regions.
We modified LAPACK's {\ttfamily{DGEQP3}} routine to obtain a truncated rank revealing QR factorization based on relative error threshold. Experiments were conducted on a system described in Table \ref{table:machine-specs}.
\begin{table}[h]
    \centering
    \caption{Details of system used for experiments}
    \label{table:machine-specs}
    \small
    \begin{tabular}{ll}
        \toprule
        & Dual AMD EPYC\textsuperscript{\texttrademark} 7502 \\
        \midrule
        Clock speed & 2.5 GHz \\
        \# cores & 2 x 32 = 64 \\
        Peak performance & 2560 GFlop/s \\
        Memory & 500 GB \\
        Compiler suite & GCC 8.4 \\
        BLAS \& LAPACK library & Intel MKL 2020.1.217 \\
        Multithreading & OpenMP 4.5 \\
        DGEMM performance & 1861 GFlop/s \\
        \bottomrule
    \end{tabular}
\end{table}

The following algorithms have been compared:
\begin{itemize}
    \item \textbf{Dense Householder}: Householder QR factorization subroutine of Intel MKL ({\ttfamily{DGEQRF}}).
    \item \textbf{Blocked MGS}: Blocked modified Gram-Schmidt-based QR factorization of weakly admissible BLR-matrices explained in Section \ref{subsec:blr_mbgs}.
    \item \textbf{Blocked Householder}: Blocked Householder-based QR factorization of BLR-matrices explained in Section \ref{subsec:blr_bch}.
    \item \textbf{Tiled Householder}: Tiled Householder-based QR factorization of BLR-matrices explained in Section \ref{subsec:blr_bh}.
\end{itemize}

We evaluate the accuracy of BLR-QR methods using two metrics: one is \textit{residual} that measures the quality of the approximate factorization; the other one is \textit{orthogonality} that measures the quality of the orthogonal factor $\tilde{Q}$. Both metrics are respectively given by
\begin{displaymath}
    \text{Res} = \frac{\|\tilde{Q}\tilde{R}-A\|_{F}}{\|A\|_{F}},\;\;
    \text{Orth} = \frac{\|\tilde{Q}^{T}\tilde{Q}-I\|_{F}}{\sqrt{n}},
\end{displaymath}
where $n$ denotes the matrix size and $\sqrt{n}$ is the Frobenius norm of order $n$ identity matrix.

\subsection{Performance on Random BLR Matrices}
\label{sec:perf_random}
First, we test our methods using randomly generated BLR matrices such that each diagonal block is a random dense matrix and each off-diagonal block is a rank-$k$ matrix obtained from the outer product of two $b \times k$ random matrices, where $b$ is the chosen BLR block size. We assume weakly admissible BLR compression with error tolerance $\epsilon=10^{-10}$.

\begin{figure}[H]
    \centering
    \includegraphics[width=0.9\textwidth]{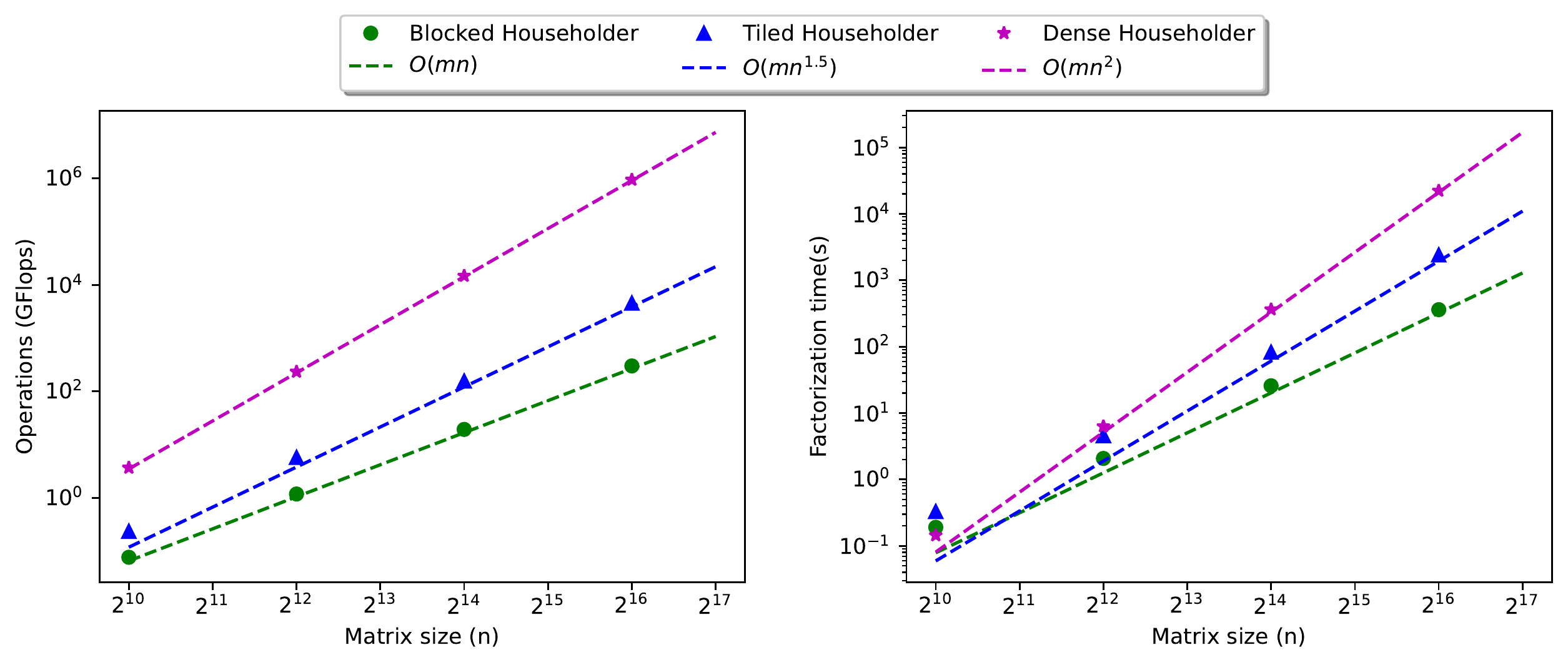}
    \caption{Flops count (left) and factorization time (right) using a single-core}
    \label{fig:res-flop-time-seq}
\end{figure}

\begin{figure}[H]
    \centering
    \includegraphics[width=0.7\textwidth]{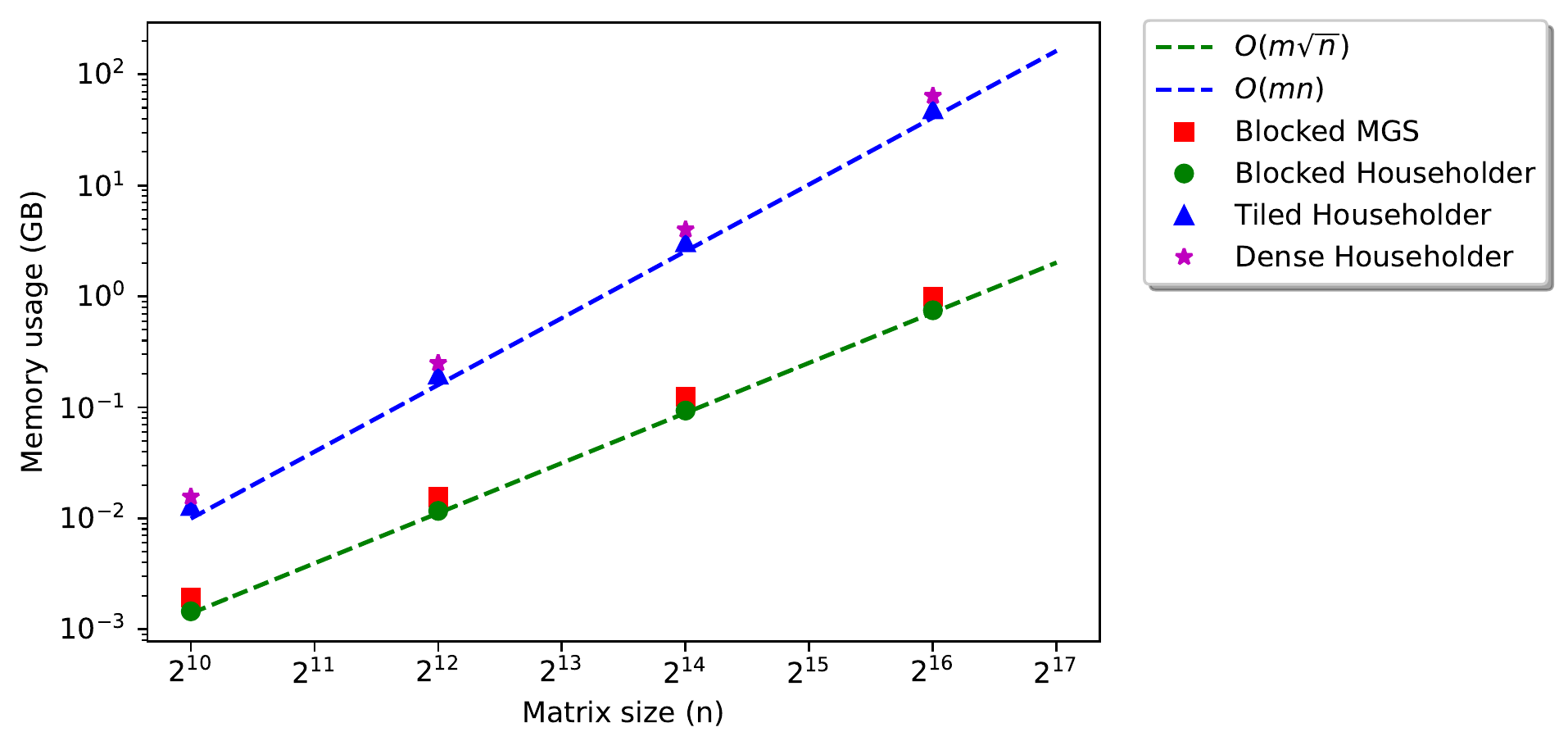}
    \caption{Memory consumption during executions using a single-core}
    \label{fig:res-memory}
\end{figure}

We first show the operation and memory complexities of our algorithms using BLR matrices of varying sizes. We generate $m \times n$ ($m=2n$) random BLR matrices with block size $b=2\sqrt{n}$ and rank $k=1$ off-diagonal blocks. Table \ref{table:res-acc-random} shows that our BLR methods produce approximate factorization accurately to the level of the prescribed error tolerance. Figure \ref{fig:res-flop-time-seq} shows the floating-point operations (flops) count and factorization time using a single-core of the machine, and Figure \ref{fig:res-memory} shows the corresponding memory consumption.
\begin{table}[H]
    \centering
    \caption{Accuracy on random BLR matrices ($\epsilon=10^{-10}$)}
    \label{table:res-acc-random}
    \small
    \begin{tabular}{rr|ll|ll}
        \toprule
        \multicolumn{1}{c}{\multirow{2}{*}{$m$}} &
        \multicolumn{1}{c}{\multirow{2}{*}{$n$}} & \multicolumn{2}{|c}{Blocked Householder} & \multicolumn{2}{|c}{Tiled Householder} \\
        & & \multicolumn{1}{|c}{Res} & \multicolumn{1}{c}{Orth} & \multicolumn{1}{|c}{Res} & \multicolumn{1}{c}{Orth} \\
        \midrule
        
        2,048 & 1,024 &
        $4.9 \cdot 10^{-15}$ & $3.7 \cdot 10^{-15}$ &
        $6.5 \cdot 10^{-14}$ & $4.1 \cdot 10^{-13}$ \\
        
        8,192 & 4,096 &
        $1.9 \cdot 10^{-14}$ & $8.0 \cdot 10^{-15}$ &
        $1.6 \cdot 10^{-13}$ & $1.7 \cdot 10^{-12}$ \\
        
        32,768 & 16,384 &
        $2.8 \cdot 10^{-14}$ & $1.7 \cdot 10^{-14}$ &
        $9.8 \cdot 10^{-14}$ & $1.1 \cdot 10^{-12}$ \\
        
        131,072 & 65,536 &
        $2.2 \cdot 10^{-13}$ & $3.7 \cdot 10^{-14}$ &
        $3.9 \cdot 10^{-13}$ & $5.3 \cdot 10^{-14}$ \\
        \bottomrule
    \end{tabular}
\end{table}

Figure \ref{fig:res-flop-time-seq} clearly shows that as the matrix size becomes larger, the flops count of both BLR-QR algorithms grow in accordance with our estimate in Section \ref{sec:blr_qr}. The right part of the figure shows that our BLR methods outperform Dense QR on large matrices. On the largest matrix ($n=65,536$), the BLR methods are more than an order of magnitude faster. However, for the smallest matrix ($n=1,024$), the BLR-QR methods, which operate on a set of small matrix blocks, suffer from the suboptimal performance of BLAS libraries on small data sizes. The benefit of using BLR factorization starts to appear when the matrix is large enough.

Figure \ref{fig:res-memory} shows that compression using BLR-matrix, followed by performing QR factorization on the compressed form leads to orders of magnitude smaller memory consumption compared to performing direct factorization on the dense matrix. The memory consumption of both blocked Householder and MGS-based methods grow as $\mathcal{O}(m\sqrt{n})$, which corresponds to the storage requirement of a rectangular BLR-matrix. Our blocked Householder-based method consumes slightly less memory compared to the existing blocked MGS-based method because it reuses the lower triangular part of $\tilde{R}$ (for $Y$ matrices) plus a block diagonal matrix (for $T$ matrices) to implicitly store the orthogonal factor $\tilde{Q}$, whereas the MGS-based method explicitly forms two BLR-matrices $\tilde{Q}$ and $\tilde{R}$ during the factorization. However, the tiled Householder-based method consumes more memory than the other BLR methods since it needs more additional space to store the $T$ matrices coming from the update QR factorization steps. This leads to a memory consumption that grows similarly to the dense QR. As a remedy, an inner blocking technique \cite{Buttari2009_ParTiledDense, Quintana2009_QR} could be employed to reduce additional storage for the $T$ matrices.

\begin{figure}[H]
    \centering
    \includegraphics[width=0.9\textwidth]{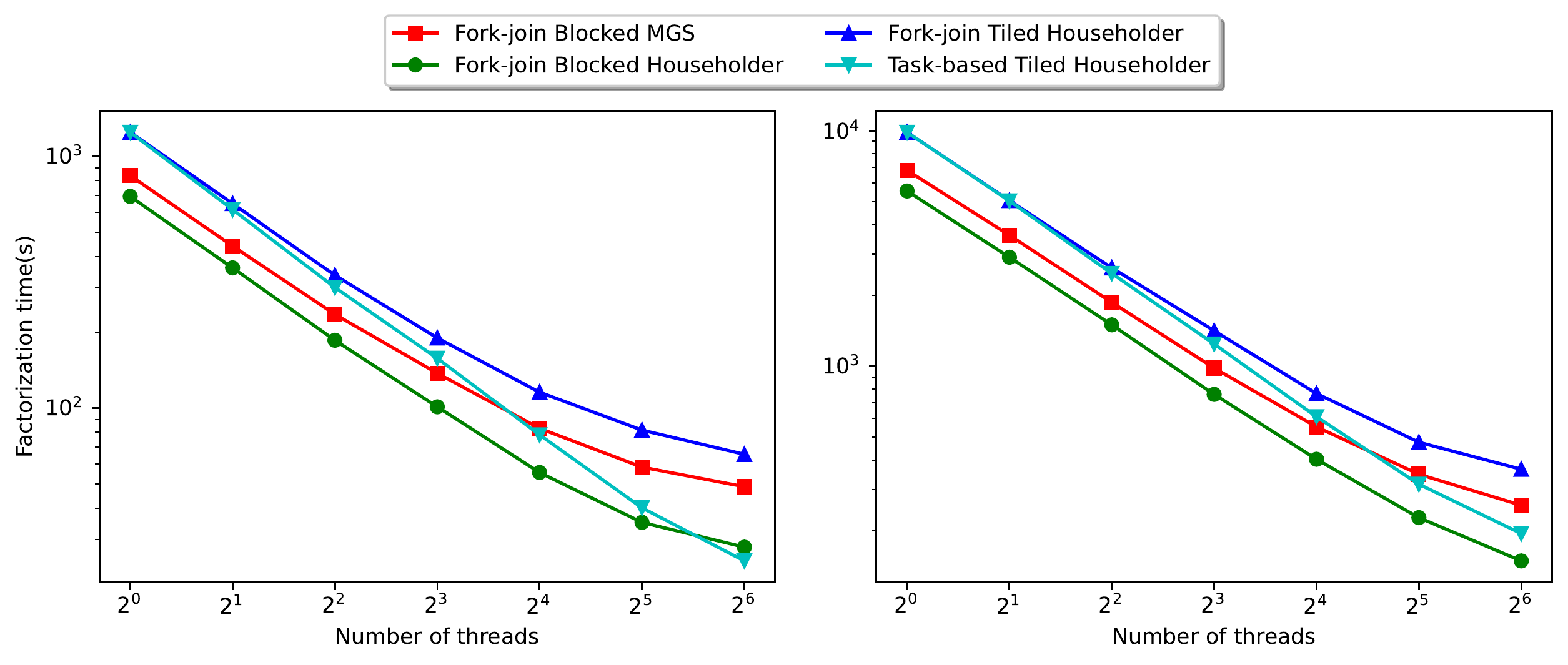}
    \caption{Factorization time using different number of threads: $n$=16,384 (left); $n$=32,768 (right)}
    \label{fig:res-weak-par-time}
\end{figure}

\begin{figure}[H]
    \centering
    \includegraphics[width=0.9\textwidth]{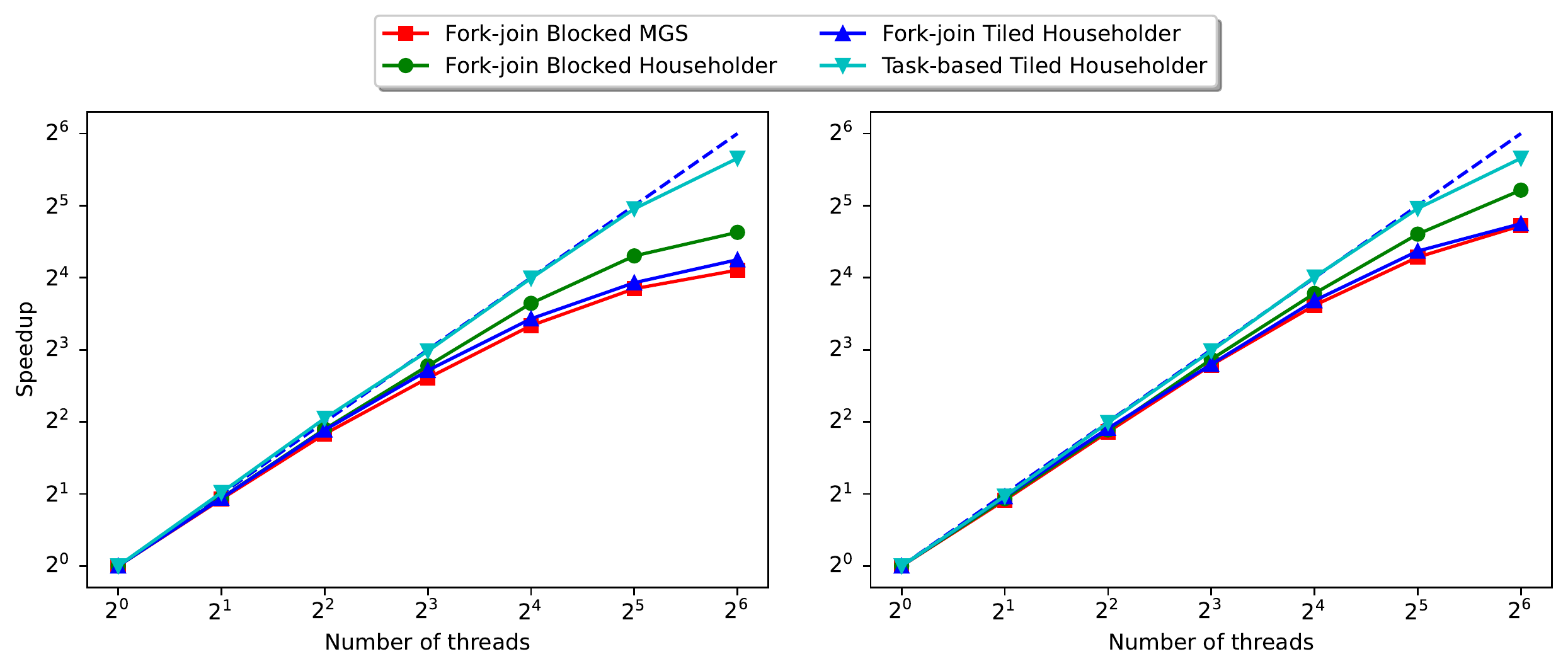}
    \caption{Speedup using different number of threads: $n$=16,384 (left); $n$=32,768 (right)}
    \label{fig:res-weak-par-speedup}
\end{figure}

\begin{figure}[h]
    \centering
    \includegraphics[width=0.9\textwidth]{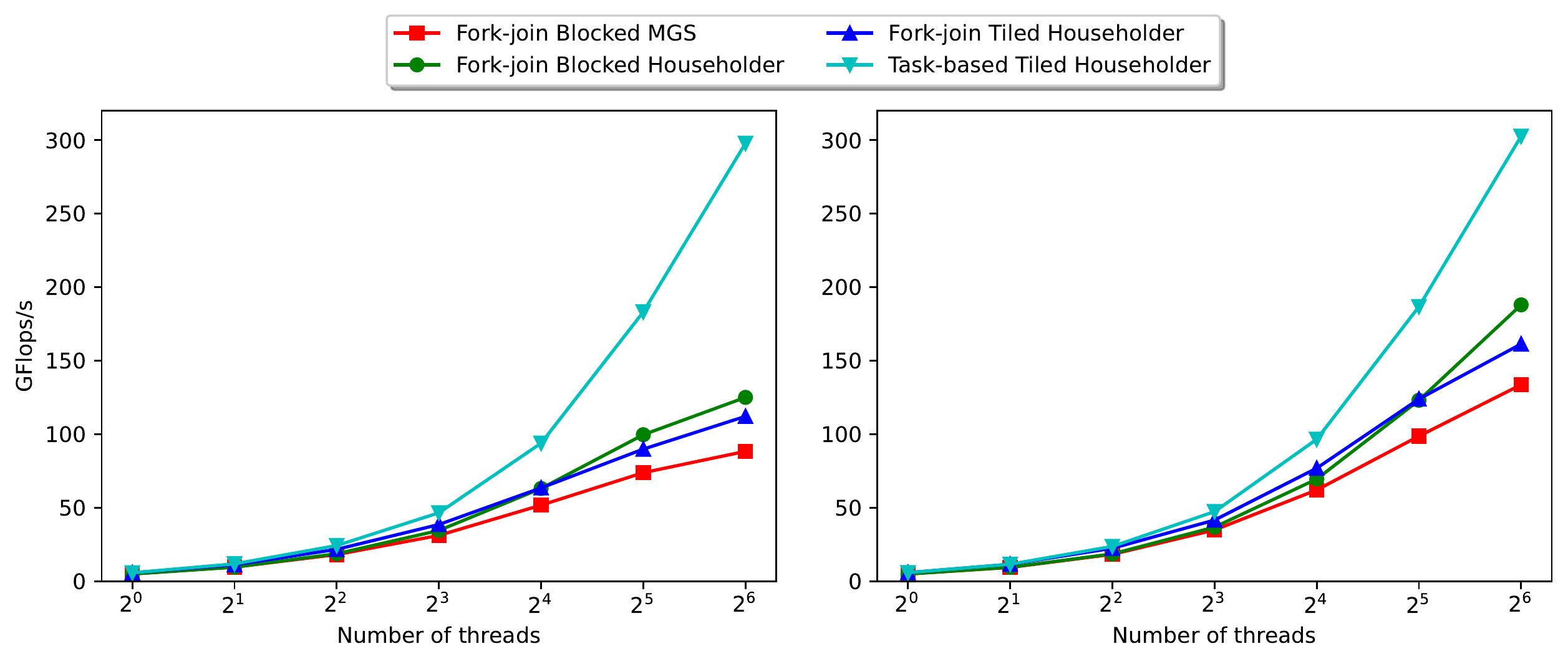}
    \caption{Flops rate using different number of threads: $n$=16,384 (left); $n$=32,768 (right)}
    \label{fig:res-weak-par-flopsrate}
\end{figure}

We now demonstrate the parallel scalability of our methods. We use two random $m \times n$ ($m=2n$) BLR-matrices with $n=16,384$ and $n=32,768$, block size $b=256$, and rank $k=16$ off-diagonal blocks. We compare our parallel algorithms with the existing parallel blocked MGS-based algorithm. Figure \ref{fig:res-weak-par-time}, \ref{fig:res-weak-par-speedup}, and \ref{fig:res-weak-par-flopsrate} show the factorization time, speedup, and flops rate using different number of threads, respectively.

Figure \ref{fig:res-weak-par-time} shows that all BLR methods scale nicely using up to 64 cores of the machine. For the matrix of size $n=16,384$, the task-based tiled Householder method outperforms the blocked Householder when using 64 cores. However, on the larger matrix ($n=32,768$), using 64 cores of our machine is not sufficient for this to happen. But we can expect that when the number of threads increases, the tiled method would eventually outperform the blocked method. This shows that the finer granularity of the tiled Householder method that allows for efficient dynamic task-based execution is able to overcome the induced extra operations once we have a large number of computing cores.

Figure \ref{fig:res-weak-par-speedup} shows that the fork-join blocked MGS, tiled Householder, and blocked Householder-based methods showed a speedup of up to 26, 26, and 37 times, respectively. The fork-join blocked Householder-based methods showed higher speedups compared to the MGS-based. 
Even though both of these methods perform similar block-column-wise QR, the blocked MGS has a bottleneck of computing the $\tilde{R}_{j,k}$ (line 3-5 of Algorithm \ref{alg:forkjoin_blr_mbgs_qr}) \cite{IdaNakashima2019}.
Furthermore, the blocked Householder outperforms tiled Householder when using fork-join execution model. On the other hand, the task-based tiled Householder achieved up to 50 times speedup, thanks to the DAG-based execution that fully utilizes the dependency between operations in the tiled Householder QR, allowing it to scale almost perfectly as the number of threads increases.

Although the scalability is promising, the actual performance of the BLR-QR algorithms is still far from the peak performance of the machine, as shown in Figure \ref{fig:res-weak-par-flopsrate}. There are two reasons behind this. First, BLR-QR methods have to deal with low-rank block operations, which involve manipulating a collection of small matrices instead of one large matrix, making it more memory-bound \cite{Rouet2016_STRUMPACK, Yu2019}. Second is the suboptimal performance of BLAS routines on small data sizes. Therefore we cannot expect these algorithms to reach the same flops rate as the traditional dense algorithms.
\begin{figure}[H]
    \centering
    \includegraphics[width=0.6\textwidth]{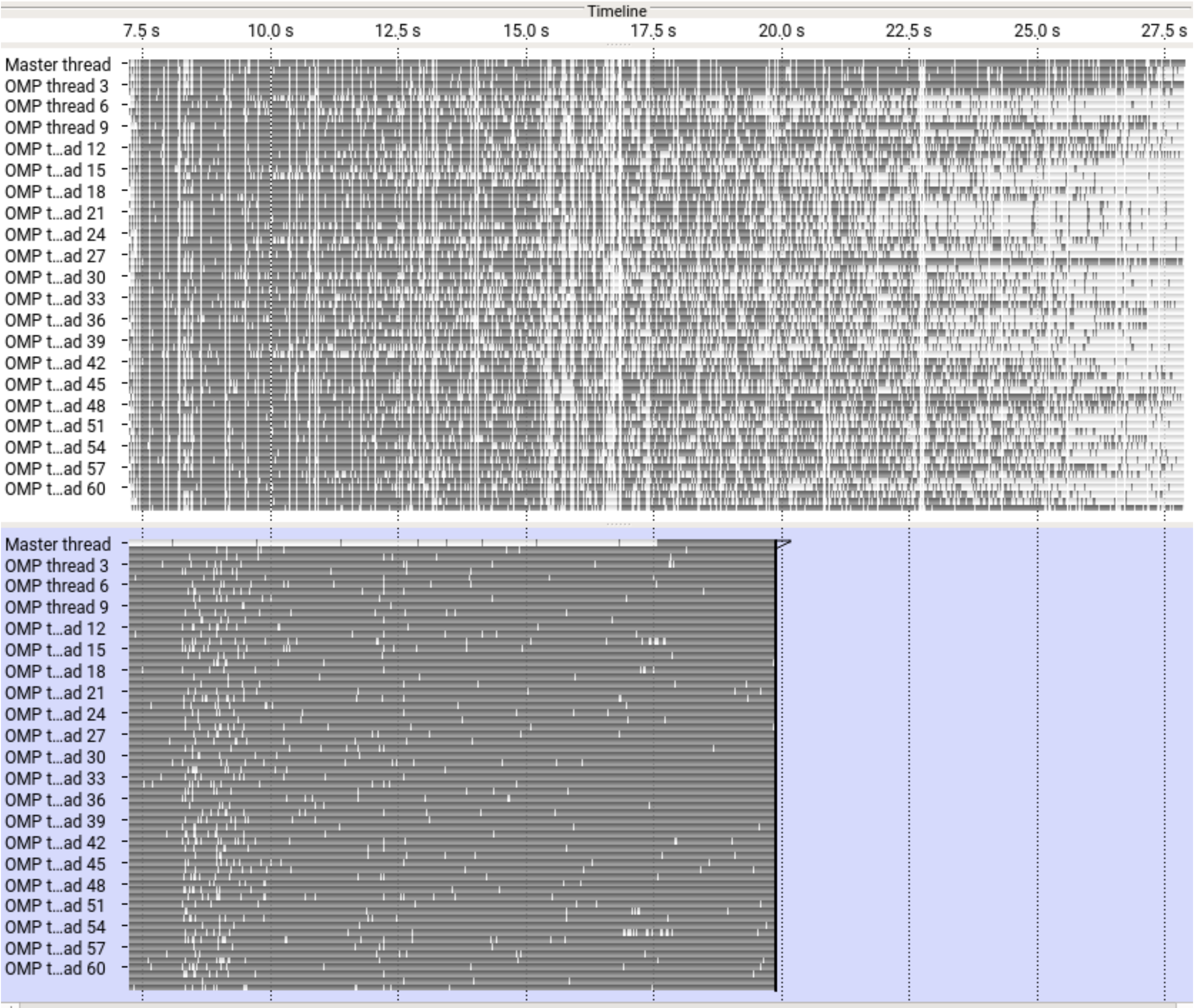}
    \caption{Execution traces of parallel tiled Householder BLR-QR on 8kx8k matrix: fork-join (top, 20.8s); task-based (bottom, 12.9s). BLR-QR executions start at around 7.5 seconds of the timeline}
    \label{fig:res-par-traces}
\end{figure}

Samples of a parallel execution trace are shown in Figure \ref{fig:res-par-traces} where the grey part corresponds to computation and the white part corresponds to synchronization and overhead. The fork-join model has many synchronizations involving all threads, which means threads that completed their task first need to wait for the others before proceeding with the execution. However, the task-based execution showed very loose synchronization because once a thread finishes a task, it can take another "ready" task from the task pool and begin another execution without waiting for other threads. This leads to an out-of-order execution that eliminates unnecessary synchronizations and significantly reduces the idle time of threads. Furthermore, one can also expect the tiled method to perform better in distributed memory systems due to its fine granularity that would lead to smaller communication overhead compared to the blocked method. Note that due to the limitation of OpenMP, in our implementation, the dependency graph is constructed on the fly by the master thread, which corresponds to the large overhead in the beginning. This is not very efficient when using a small number of threads since the master thread spends more than 70\% of its time generating tasks. This overhead however is not significant when using a large number of threads.
\subsection{Accuracy on Ill-Conditioned Matrices}
\label{sec:perf_bem2d}
In this example, we demonstrate the numerical stability of our methods in factorizing ill-conditioned matrices. We use matrices arising from Boundary-Element-Method discretization of Single-Layer Potential (SLP) operator on the unit circle, generated using {\ttfamily{H2Lib}} \cite{H2Lib}. The resulting square matrices are ill-conditioned and have off-diagonal blocks with small ranks, hence we assume weakly admissible BLR compression. We set the block size $b=2\sqrt{n}$ and error tolerance $\epsilon=10^{-9}$. 

We compare the accuracy of our methods with the existing blocked MGS-based method. Table \ref{table:res-acc-bem2d} shows that as the condition number increases, our Householder methods are robust to this increase and produce numerical orthogonality on the level of the prescribed tolerance. On the other side, the orthogonality produced by the MGS method clearly deteriorates as the condition number increases.
\begin{table}[H]
    \centering
    \caption{Accuracy on ill-conditioned matrices ($\epsilon=10^{-9}$)}
    \label{table:res-acc-bem2d}
    \small
    \begin{tabular}{rrcc|ll|ll|ll}
        \toprule
        \multicolumn{1}{c}{\multirow{2}{*}{$n$}} & \multicolumn{1}{c}{\multirow{2}{*}{$\kappa_F(A)$}} & \multicolumn{1}{c}{Block} & \multicolumn{1}{c|}{Max} & \multicolumn{2}{c|}{Blocked Householder} & \multicolumn{2}{c|}{Tiled Householder} & \multicolumn{2}{c}{Blocked MGS} \\
        
        & & \multicolumn{1}{c}{Size} & \multicolumn{1}{c|}{Rank} & \multicolumn{1}{c}{Res} & \multicolumn{1}{c|}{Orth} &
        \multicolumn{1}{c}{Res} & \multicolumn{1}{c|}{Orth} & \multicolumn{1}{c}{Res} & \multicolumn{1}{c}{Orth} \\
        \midrule
        1,024 & $2.8 \cdot 10^{5}$ & 64 & 11 & 
        $6.8 \cdot 10^{-10}$ & $6.9 \cdot 10^{-11}$ &
        $6.1 \cdot 10^{-10}$ & $5.2 \cdot 10^{-11}$ &
        $5.1 \cdot 10^{-10}$ & $1.9 \cdot 10^{-8}$ \\
        
        4,096 & $4.6 \cdot 10^{6}$ & 128 & 12 & 
        $1.0 \cdot 10^{-9}$ & $1.2 \cdot 10^{-10}$ & 
        $9.6 \cdot 10^{-10}$ & $4.5 \cdot 10^{-11}$ &
        $8.6 \cdot 10^{-10}$ & $1.0 \cdot 10^{-7}$ \\
        
        16,384 & $7.4 \cdot 10^{7}$ & 256 & 12 & 
        $2.1 \cdot 10^{-9}$ & $6.2 \cdot 10^{-11}$ &
        $1.8 \cdot 10^{-9}$ & $6.7 \cdot 10^{-11}$ &
        $1.8 \cdot 10^{-9}$ & $1.5 \cdot 10^{-6}$ \\
        
        32,768 & $2.9 \cdot 10^{8}$ & 512 & 13 & 
        $2.3 \cdot 10^{-9}$ & $6.0 \cdot 10^{-11}$ & 
        $2.0 \cdot 10^{-9}$ & $5.3 \cdot 10^{-11}$ &
        $2.0 \cdot 10^{-9}$ & $6.9 \cdot 10^{-6}$ \\
        \bottomrule
    \end{tabular}
\end{table}
\subsection{Performance on Spatial Statistics Problems}
\label{sec:perf_gaussian3d}
In this example, we use square matrices arising from the Spatial Statistics problem with exponential kernel on uniform 3D grids, generated using {\ttfamily{STARS-H}} \cite{Stars-H}. The resulting matrices have many off-diagonal blocks with relatively high ranks, thus a strongly admissible compression is preferred.

\begin{figure}[H]
    \centering
    \includegraphics[width=0.9\textwidth]{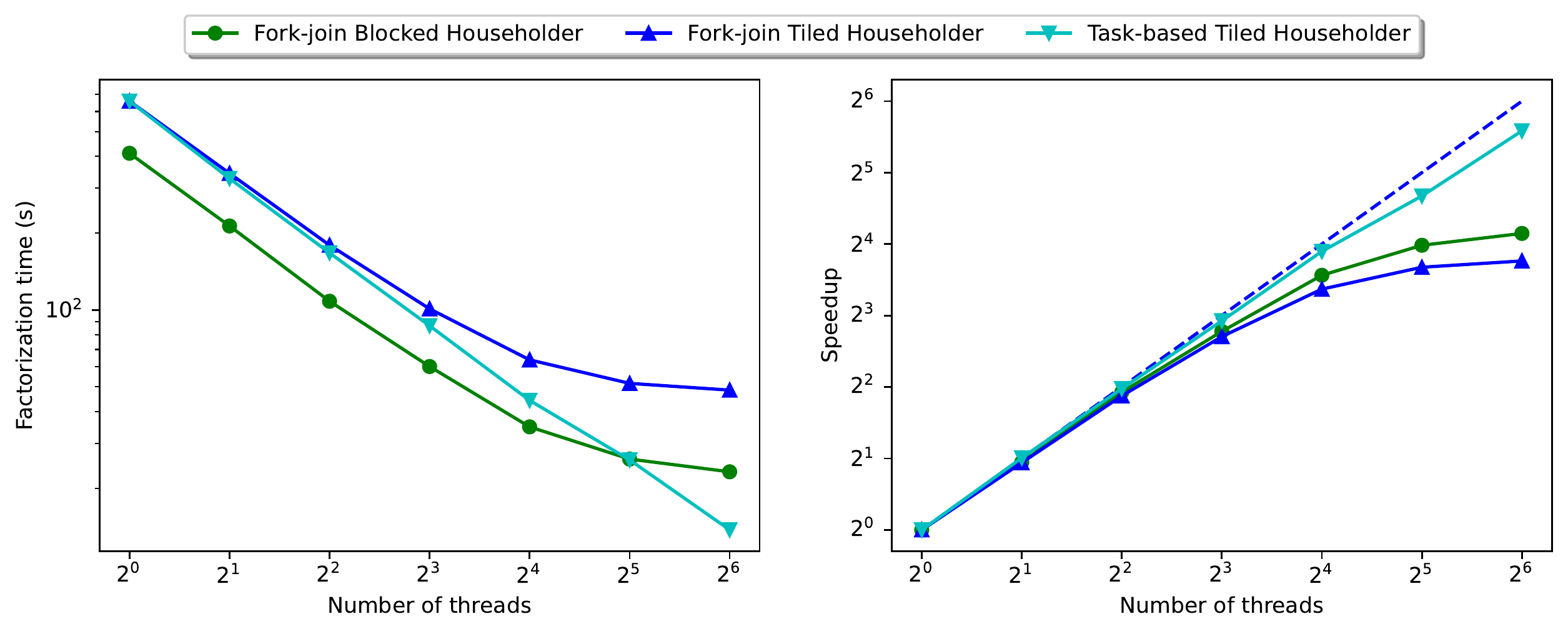}
    \caption{Parallel scalability on spatial statistics problem: factorization time (left); speedup (right)}
    \label{fig:res-gaussian3d-par-time-speedup}
\end{figure}

We first demonstrate the parallel scalability of our methods for the matrix of order $n=16,384$ compressed with block size $b=256$, tolerance $\epsilon=10^{-6}$, and admissibility constant $\eta=0.3$. Figure \ref{fig:res-gaussian3d-par-time-speedup} shows that using 64 cores the fork-join tiled and blocked householder methods achieve a speedup of 13 and 17 times, respectively. On the other hand, the task-based tiled Householder method achieves 47 times speedup, allowing it to become faster than the blocked method.

Next, we evaluate the performance and accuracy of our methods using matrices of varying sizes and admissibility constant. We compare our fork-join blocked and task-based tiled Householder with the parallel dense QR of Intel MKL using 64 cores in terms of factorization time and memory consumption. Table \ref{table:res-acc-exp3d} shows that our BLR-QR methods are able to produce residual and orthogonality to the level of the prescribed error tolerance. It also shows that we can fine-tune the admissibility constant (introduce more/less off-diagonal inadmissible block, i.e. decrease/increase the maximum rank of admissible off-diagonal blocks) to reach the optimal factorization time and memory consumption. On the matrix of order 16k, the BLR methods lose to the dense QR. However, for the larger matrix of order 65k, the BLR methods are faster and achieve up to 80\% less memory usage compared to the dense QR.
\begin{table}[H]
    \centering
    \caption{Accuracy on 3D Spatial Statistics problems ($\epsilon=10^{-6}$)}
    \label{table:res-acc-exp3d}
    \small
    \resizebox{\textwidth}{!}{%
    \begin{tabular}{r|rr|ccc|rrll|rrll}
        \toprule
        & \multicolumn{2}{|c}{Dense QR} & \multicolumn{3}{|c}{BLR} & \multicolumn{4}{|c}{Blocked Householder} & \multicolumn{4}{|c}{Tiled Householder} \\

        \multicolumn{1}{c}{\multirow{2}{*}{$n$}} & \multicolumn{1}{|c}{Mem} & \multicolumn{1}{c}{Time} & \multicolumn{1}{|c}{Block} & \multicolumn{1}{c}{$\eta$} & \multicolumn{1}{c}{Max} & \multicolumn{1}{|c}{Mem} & \multicolumn{1}{c}{Time} & \multicolumn{1}{c}{Res} & \multicolumn{1}{c}{Orth} & \multicolumn{1}{|c}{Mem} & \multicolumn{1}{c}{Time} & \multicolumn{1}{c}{Res} & \multicolumn{1}{c}{Orth}\\
        
        & \multicolumn{1}{|c}{(MB)} & \multicolumn{1}{c}{(s)} & \multicolumn{1}{|c}{Size} & & \multicolumn{1}{c}{Rank} & \multicolumn{1}{|c}{(MB)} & \multicolumn{1}{c}{(s)} & & & \multicolumn{1}{|c}{(MB)} & \multicolumn{1}{c}{(s)} & & \\
        \midrule

        \multirow{3}{*}{16,384} & \multirow{3}{*}{2,048} & \multirow{3}{*}{8.26} & 256 & 0.2 & 117 & 847 & 32.5 & $1.1 \cdot 10^{-9}$ & $1.3 \cdot 10^{-11}$ & 1,857 & 16.9 & $1.5 \cdot 10^{-9}$ & $1.1 \cdot 10^{-9}$ \\
        & & & 256 & 0.3 & 86 & 954 & 23 & $8.0 \cdot 10^{-10}$ & $1.1 \cdot 10^{-11}$ & 1,964 & 13.7 & $1.4 \cdot 10^{-9}$ & $1.0 \cdot 10^{-9}$ \\
        & & & 256 & 0.4 & 58 & 1,284 & 16.4 & $2.1 \cdot 10^{-10}$ & $3.8 \cdot 10^{-12}$ & 2,293 & 10.8 & $3.8 \cdot 10^{-10}$ & $2.2 \cdot 10^{-10}$ \\
        \midrule
        
        \multirow{3}{*}{65,536} & \multirow{3}{*}{32,768} & \multirow{3}{*}{310} & 512 & 0.2 & 175 & 6,698 & 257.7 & $3.7 \cdot 10^{-9}$ & $1.5 \cdot 10^{-10}$ & 22,957 & 215.9 & $6.3 \cdot 10^{-9}$ & $4.0 \cdot 10^{-9}$ \\
        & & & 512 & 0.3 & 94 & 10,236 & 185.6 & $1.5 \cdot 10^{-9}$ & $8.8 \cdot 10^{-11}$ & 26,491 & 218.8 & $2.8 \cdot 10^{-9}$ & $1.8 \cdot 10^{-9}$ \\
        & & & 512 & 0.4 & 51 & 14,684 & 230.9 & $5.1 \cdot 10^{-10}$ & $4.3 \cdot 10^{-11}$ & 30,937 & 286.9 & $9.9 \cdot 10^{-10}$ & $6.3 \cdot 10^{-10}$ \\
        \bottomrule
    \end{tabular}%
    }
\end{table}
\subsection{Performance on Inverse Poisson Problems}
\label{sec:perf_invpoi2d}
In this example, we use sparse least squares matrices arising from the Inverse Poisson problem defined on uniform 2D grids, generated using the MATLAB code of {\ttfamily{spaQR}} \cite{Gnanasekaran2021_spaQR_LS}. We use strongly admissible BLR compression with block size $b=2\sqrt{n}$ and tolerance $\epsilon=10^{-10}$. Since the geometry information is not available, we attempt to compress every off-diagonal block and revert back to dense the blocks whose rank is larger than $b/2$.

\begin{figure}[h]
    \centering
    \includegraphics[width=0.9\textwidth]{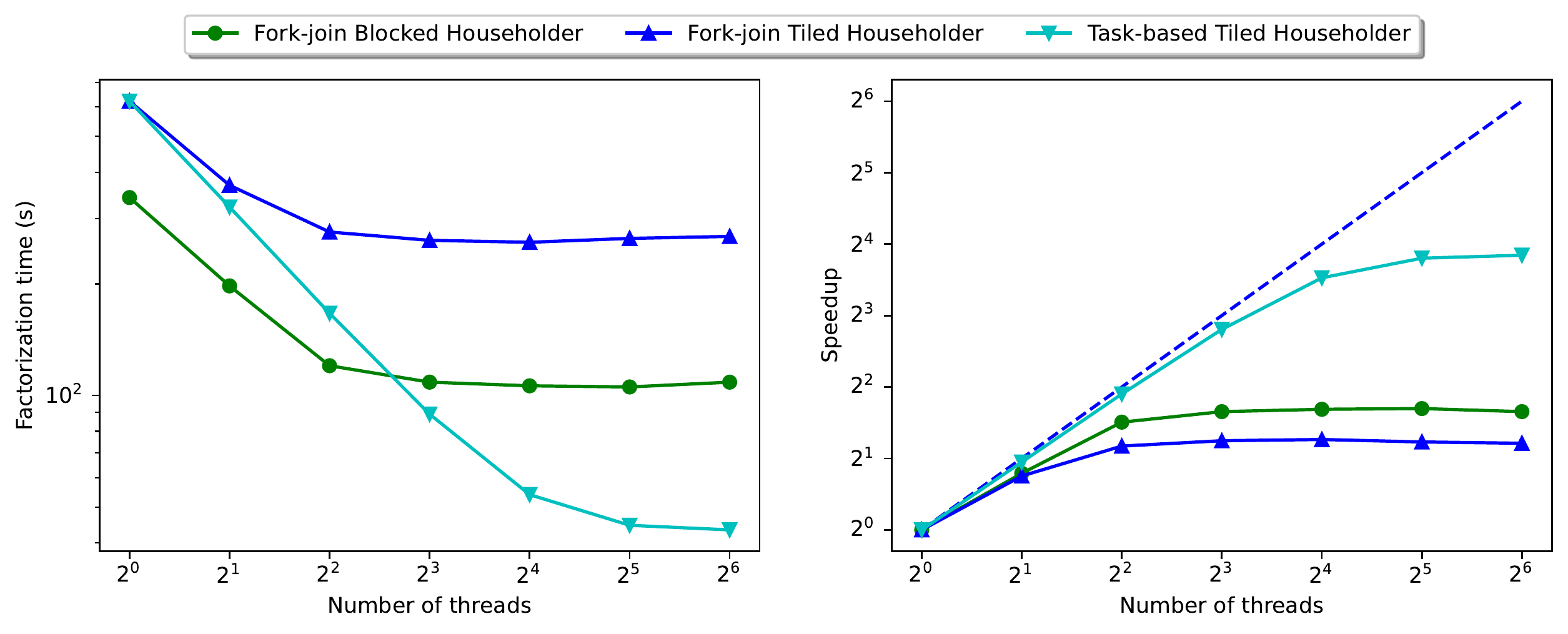}
    \caption{Parallel scalability on sparse least squares problem: factorization time (left); speedup (right)}
    \label{fig:res-invpoi2d-par-time-speedup}
\end{figure}

Figure \ref{fig:res-invpoi2d-par-time-speedup} shows the parallel scalability of our BLR methods for the matrix of size $74,112 \times 36,864$. Using 64 cores, the fork-join tiled, fork-join blocked, and task-based tiled Householder methods achieve a speedup of 2, 3, and 14 times, respectively. These relatively lower speedups come from the fact that although the dimension is quite large, the resulting BLR matrix is dominated by zero blocks that make the actual computation load smaller.

We then compare our fork-join blocked and task-based tiled Householder with the parallel dense QR of Intel MKL using all 64 cores of the machine. Table \ref{table:res-acc-invpoi2d} shows that the BLR methods are faster than Dense QR while producing residual and orthogonality to the level of the prescribed error tolerance. Moreover, the tiled method is faster than the blocked method in all three sparse matrices that we use.
\begin{table}[H]
    \centering
    \caption{Accuracy on 2D Inverse Poisson problems ($\epsilon=10^{-10}$)}
    \label{table:res-acc-invpoi2d}
    \small
    \resizebox{\textwidth}{!}{%
    \begin{tabular}{rr|rr|rrll|rrll}
        \toprule
        & & \multicolumn{2}{|c}{Dense QR} & \multicolumn{4}{|c}{Blocked Householder} & \multicolumn{4}{|c}{Tiled Householder} \\
        
        \multicolumn{1}{c}{\multirow{2}{*}{$m$}} & \multicolumn{1}{c}{\multirow{2}{*}{$n$}} & \multicolumn{1}{|c}{Mem} & \multicolumn{1}{c}{Time} & \multicolumn{1}{|c}{Mem} & \multicolumn{1}{c}{Time} & \multicolumn{1}{c}{Res} & \multicolumn{1}{c}{Orth} & \multicolumn{1}{|c}{Mem} & \multicolumn{1}{c}{Time} & \multicolumn{1}{c}{Res} & \multicolumn{1}{c}{Orth} \\
        
        & & \multicolumn{1}{|c}{(MB)} & \multicolumn{1}{c}{(s)} & \multicolumn{1}{|c}{(MB)} & \multicolumn{1}{c}{(s)} & & & \multicolumn{1}{|c}{(MB)} & \multicolumn{1}{c}{(s)} & & \\
        \midrule
        74,112 & 36,864 & 20,844 & 200 & 3,727 & 106 & $1.4 \cdot 10^{-11}$ & $1.3 \cdot 10^{-11}$ & 20,017 & 46 & $8.4 \cdot 10^{-12}$ & $5.7 \cdot 10^{-12}$ \\
        100,800 & 50,176 & 38,587 & 479 & 6,117 & 190 & $1.5 \cdot 10^{-11}$ & $1.4 \cdot 10^{-11}$ & 36,125 & 82 & $8.4 \cdot 10^{-12}$ & $5.4 \cdot 10^{-12}$ \\
        131,584 & 65,536 & 65,792 & 1569 & 9,399 & 354 & $1.6 \cdot 10^{-11}$ & $1.6 \cdot 10^{-11}$ & 60,351 & 161 & $8.6 \cdot 10^{-12}$ & $5.4 \cdot 10^{-12}$ \\
        \bottomrule
    \end{tabular}%
    }
\end{table}
\section{Conclusion}
\label{sec:conclusion}

We have presented two new algorithms for Householder QR factorization of Block Low-Rank matrices. One that performs block-column-wise QR based on the blocked Householder method, and another one that performs fine-grained, block-wise QR based on the tiled Householder method. We have shown that both algorithms exploit BLR structure to achieve arithmetic complexity of $\mathcal{O}(mn)$ and $\mathcal{O}(mn^{1.5})$, respectively. We have compared our algorithms with an existing BLR-QR method that is based on the blocked Modified Gram Schmidt iteration. We also compared them to a state-of-the-art vendor-optimized dense Householder QR of Intel MKL. Numerical experiments showed that all BLR methods are more than an order of magnitude faster than the dense QR of MKL. The BLR methods are also more efficient in terms of memory consumption, possibly saving hundreds of gigabytes of memory for huge matrices.

We also have demonstrated the parallelization of our algorithms using both traditional fork-join and modern task-based execution models. We compared our parallel algorithms with an existing fork-join blocked MGS-based parallel algorithm. Results showed that our task-based tiled Householder algorithm outperforms the fork-join methods, thanks to the dynamic task-based execution that allows for out-of-order execution with very loose synchronization between the threads. We have shown that in a shared-memory parallel environment, the benefit that comes from a finely grained algorithm is able to overcome the extra operations that it introduces.

Numerical experiments also showed that our methods can be used in various computational science problems to produce approximate QR factorization with controllable accuracy. This shows that BLR matrices can provide a good approximation for the resulting orthogonal and upper triangular factors.
Both Householder and MGS-based BLR methods produced approximate QR factorization of similar residual. However in terms of orthogonality, our method is robust to ill-conditioning, whereas the existing MGS-based method suffers from numerical instability.
\section*{Acknowledgements}
This work was supported by JSPS KAKENHI Grant Number JP20K20624 and JP21H03447.
This work is supported by "Joint Usage/Research Center for Interdisciplinary Large-scale Information Infrastructures" and "High Performance Computing Infrastructure" in Japan (Project ID: jh210024-NAHI).
This work is conducted as research activities of AIST - Tokyo Tech Real World Big-Data Computation Open Innovation Laboratory (RWBC-OIL).

\bibliographystyle{plain}
\bibliography{0_0_main}

\end{document}